\documentclass[11pt]{article}
\usepackage{amsmath}
\usepackage{amssymb,amsmath}
\usepackage{mathrsfs}
\usepackage{graphicx}
\usepackage{color}
\let\f=\frac
\let\p=\psi

\let\wt=\widetilde
\let\wh=\widehat

\def\cC{{\cal C}}

\def\cF{{\cal F}}

\def\cP{{\cal P}}
\def\cQ{{\cal Q}}

\def\cS{{\cal S}}

\def\cZ{{\cal Z}}

\def\wh{\widehat}

\def\p{\partial}

\def\wt{\widetilde}
\def\ov{\overline}

\def\eqdefa{\buildrel\hbox{\footnotesize def}\over =}

\def\C{\mathop{\bf C\kern 0pt}\nolimits}
\def\DD{\mathop{\bf D\kern 0pt}\nolimits}
\def\K{\mathop{\bf K\kern 0pt}\nolimits}
\def\N{\mathop{\bf N\kern 0pt}\nolimits}
\def\Q{\mathop{\bf Q\kern 0pt}\nolimits}
\def\R{\mathop{\bf R\kern 0pt}\nolimits}

\def\ddq{\dot \Delta_q}

\renewcommand{\div}{\mbox{\rm div}\;\!}

\newcommand{\beq}{\begin{equation}}
\newcommand{\eeq}{\end{equation}}
\newcommand{\ben}{\begin{eqnarray}}
\newcommand{\een}{\end{eqnarray}}
\newcommand{\beno}{\begin{eqnarray*}}
\newcommand{\eeno}{\end{eqnarray*}}

\newtheorem{Theorem}{Theorem}[section]
\newtheorem{Definition}[Theorem]{Definition}
\newtheorem{Proposition}[Theorem]{Proposition}
\newtheorem{Lemma}[Theorem]{Lemma}
\newtheorem{Corollary}[Theorem]{Corollary}
\newtheorem{Remark}[Theorem]{Remark}

\numberwithin{equation}{section}

\allowdisplaybreaks \numberwithin{equation} {section}

\begin{document}
\title{The unique global solvability  and  optimal time  decay rates for a  multi-dimensional  compressible generic two-fluid
model with capillarity effects  \thanks {Research supported by the
National Natural Science Foundation of China (11501332,11771043,11371221), the  Natural Science Foundation of Shandong Province (ZR2015AL007),
 and Young Scholars Research Fund of Shandong University of Technology.}
}
\author{ Fuyi  Xu$^\dag$   \  Meiling  Chi\ \\[2mm]
 { \small  School of Mathematics and  Statistics, Shandong University of Technology,}\\
  { \small Zibo,    255049,  Shandong,    China}}
         \date{}
         \maketitle
\noindent{\bf Abstract}\ \ \ The present paper deals with
the Cauchy problem of a compressible  generic two-fluid
model with capillarity effects in any dimension $N\geq2$. We first study the  unique global solvability  of the model  in spaces with critical regularity indices with respect to the scaling of the associated equations.
Due to the presence of the  capillary terms, we exploit the parabolic properties of the linearized system for all frequencies which enables us to apply  contraction mapping principle to show the unique global solvability  of  strong
solutions  close to a stable
equilibrium state. Furthermore,  under a mild additional decay assumption involving only the low frequencies of the  data,  we establish  the optimal  time decay rates  for  the  constructed global
solutions.
\vskip   0.2cm \footnotetext[1]{$^\dag$Corresponding author.}
\vskip   0.2cm \footnotetext[2]{E-mail addresses: zbxufuyi@163.com(F.Xu), \ chimeiling0@163.com(M. Chi).} \setlength{\baselineskip}{20pt}

\section{Introduction and Main Results}
\setcounter{section}{1}\setcounter{equation}{0} \ \ \ \ \
 In the present paper, we consider  the following   multi-dimensional $(N\geq2)$ non-conservative viscous
compressible two-fluid model with capillarity effects:
\begin{align}\label{equ:CTFS}
\left\{
\begin{aligned}
&\alpha^{+}+\alpha^{-}=1,\\
&\p_t(\alpha^{\pm}\rho^{\pm})+\textrm{div}(\alpha^{\pm}\rho^{\pm}u^{\pm})=0, \\
&\p_t(\alpha^{\pm}\rho^{\pm}u^{\pm})+\textrm{div}(\alpha^{\pm}\rho^{\pm}
u^{\pm}\otimes u^{\pm})+\alpha^{\pm}\nabla P^{\pm}(\rho^{\pm})-\sigma^{\pm}\alpha^{\pm}\rho^{\pm}\nabla \Delta(\alpha^{\pm}\rho^{\pm})
=\textrm{div}(\alpha^{\pm}\tau^{\pm}), \\
&P^{+}(\rho^{+})=A^{+}(\rho^{+})^{\bar{\gamma}^{+}}=P^{-}(\rho^{-})=A^{-}(\rho^{-})^{\bar{\gamma}^{-}},
\end{aligned}
\right.
\end{align}
where the variable $0\leq\alpha^{+}(x,t)\leq1$ is the volume fraction of fluid $+$ in one of the two gases, and $0\leq\alpha^{-}(x,t)\leq1$  is the  volume fraction of the other fluid$-$. Moreover, $\rho^{\pm}(x,t)\geq0$, $u^{\pm}(x,t)$ and $P^{\pm}(\rho^{\pm})=A^{\pm}(\rho^{\pm})^{\overline{\gamma}^{\pm}}$ are, respectively, the
densities, the velocities, and the two pressure functions of the fluids. $\sigma^{\pm}>0$ is the capillary coefficients. It is assumed
that $\overline{\gamma}^{\pm}>1,A^{\pm}>0$ are constants. In what follows, we set $A^{+}=A^{-}=1$ without loss of any generality. Also, $\tau^{\pm}$ are the viscous stress tensors
\begin{equation}\label{1.2}
\tau^{\pm}:=2\mu^{\pm}D(u^{\pm})+\lambda^{\pm}\textrm{div}u^{\pm}\textrm{Id},
\end{equation}
where $D(u^{\pm})\eqdefa\frac{\nabla u^{\pm}+\nabla^{t}u^{\pm}}{2}$ stand for the deformation tensor,  the constants  $\mu^{\pm}$ and $\lambda^{\pm}$ are the (given) shear and bulk viscosity coefficients satisfying $\mu^{\pm}>0$ and $\lambda^{\pm}+2\mu^{\pm}>0$.   This system is known as a two-fluid flow system with algebraic closure,  which is widely used in industrial applications, such as nuclear, power, oil-and-gas, micro-technology and so on,  and  we refer readers to Refs \cite{BDG,BHL,GN,EWW,NP,MT,NK1,NK2} for more discussions about this model and related models.

The system \eqref{equ:CTFS} is a highly nonlinear  partial differential   system  with
 the mixed hyperbolic-parabolic property. As a
matter of fact, there is no diffusion on the mass conservation equations, whereas velocity evolves according to the parabolic equations due to the  viscosity phenomena.
We should point out that the system \eqref{equ:CTFS} includes
 important single phase flow models such as  the compressible Navier-Stokes equations (i.e., $\sigma^{\pm}=0$) and the compressible Navier-Stokes-Korteweg model
 when one of the two phases volume fraction tends to zero (i.e., $\alpha^{+}=0$ or $\alpha^{-}=0$ ).  As  the
extremely important models to describe compressible fluids, these two systems  have attracted a lot of attention among many analysts and many important theories have been developed.  Here, we briefly review some  of the most relevant papers
about global well-posedness and  large time behaviors of the  solutions for two systems. For  the compressible Navier-Stokes equations,  Lions \cite{Lions} proved the global existence of weak solutions
for large initial data. However, the question of uniqueness of weak solutions remains open,
even in the two dimensional case.  Nash \cite{NA} considered the local well-posedness for smooth data away from a vacuum. Matsumura-Nishida \cite{MN2}  first studied the global existence of smooth
solutions close to a stable
equilibrium  for the initial data $(\rho_{0}-\bar{\rho}, u_{0})\in H^{3}(\mathbb{R}^{3})\times H^{3}(\mathbb{R}^{3})$, where $\bar{\rho}$ is a positive constant.
Matsumura-Nishida \cite{MN1},  provided that the initial perturbation $(\rho_{0}-\bar{\rho}, u_{0})$ is sufficiently small in $H^{4}(\mathbb{R}^{3})\cap L^{1}(\mathbb{R}^{3})$, obtained the following the optimal time decay rates
\begin{equation}\label{eq:decate1}\|(\rho-\bar{\rho}, u)\|_{L^{2}}\leq C(1+t)^{-\frac{3}{4}},\, t\geq0.\end{equation}
Later, if the small initial disturbance belongs to $H^{s}(\mathbb{R}^{N})\cap W^{s,1}(\mathbb{R}^{N})$ with the integers $s\geq[\frac{N}{2}]+3$ and the space dimensions $N=2,3$, Ponce \cite{Ponce} established  the optimal $L^{p}$-time decay rates
\begin{equation}\label{eq:decate2}\|\nabla^{k}(\rho-\bar{\rho}, u)\|_{L^{p}}\leq C(1+t)^{-\frac{N}{2}(1-\frac{1}{p})-\frac{k}{2}},\, t\geq0,\, 0\leq k\leq 2,\, 2\leq p\leq\infty .\end{equation}
In 2007, Duan et al.\cite{Dun2} proved  the optimal time decay rates for the compressible Navier-Stokes equations
with potential forces  if the small initial disturbance belongs to  $H^{3}(\mathbb{R}^{3})$ and  the initial
perturbation is bounded in $L^{1}(\mathbb{R}^{3})$
without the smallness of $L^{1}$-norm of the initial disturbance. Note that $L^{1}(\mathbb{R}^{3})$ is included in $\dot{B}^{1}_{1,\infty}(\mathbb{R}^{3})$ (which will be defined in Sec. 2),  Li and Zhang \cite{LZ} extended the known results from \cite{MN1,MN2} and
showed the optimal time decay rates of the the density and momentum  when initial perturbation is sufficiently small in $ H^{l}(\mathbb{R}^{3})\cap\dot{B}_{1,\infty}^{-s}(\mathbb{R}^{3}), l\geq4$, and $s\in[0,1]$
\begin{equation}\label{eq:decate3}
\|\partial_{x}^{k}(\rho-\bar{\rho})(t)\|_{L^{2}}+
\|\partial_{x}^{k}m(t)\|_{L^{2}}\lesssim(1+t)^{-\frac{3}{4}-\frac{s+|k|}{2}}, \, \hbox{for}\, |k|=0,1.
\end{equation}
Recently, the present author and Chi \cite{XC} removed the smallness  condition of $\dot{B}_{1,\infty}^{-s}(\mathbb{R}^{3})$ for  $s\in[0,\frac{1}{2})$.
Concerning the optimal time decay rates of global solutions to the compressible Navier-
Stokes equations (with or without external potential force), we refer the reader to the
papers \cite{Dun1,GW,Ka1,Ka2,SK} and references therein.
In  critical framework,  Danchin \cite{Dan1} first proved   the existence  and uniqueness of the global strong solutions for  initial data  close to a stable
equilibrium state. The optimal time decay estimates issue for the compressible Navier-Stokes equations  in the critical regularity framework for  dimension $N\geq3$  has been addressed only very recently
by Okita in \cite{O},  provided the data are additionally in some
superspace of $L^1$.  In the survey paper \cite{Dan2}, Danchin proposed  another description of the
time decay which allows to handle  dimension $N\geq2$ in the $L^2$ critical framework. For  the compressible Navier-Stokes-Korteweg model,  Hattori and  Li \cite{HL1,HL2} established the local existence of smooth solutions  with large initial data
and global existence of smooth solutions around constant states
for small  initial data  $(\rho_0,u_0)$
 in  Sobolev spaces $H^{s}(\mathbb{R}^{N})\times H^{s-1}(\mathbb{R}^{N})$ with $s\geq\frac{N}{2}+4$ and $N=2,3$.  Recently, researchers in \cite{Tan2} proved the global existence and obtain the following optimal decay rates of strong solutions for small
initial data in some Sobolev spaces which have lower regularity than that of \cite{HL2} in three
dimensional case
\begin{equation}\label{eq:decate4}\|(\rho-\bar{\rho}, u)\|_{L^{2}}\leq C(1+t)^{-\frac{3}{4}},\, t\geq0,\end{equation}
\begin{equation}\label{eq:decate5}\|\nabla(\rho-\bar{\rho}, u)\|_{L^{2}}\leq C(1+t)^{-\frac{5}{4}},\, t\geq0,\end{equation}
\begin{equation}\label{eq:decate51}\|(\rho-\bar{\rho}, u)\|_{L^{p}}\leq C(1+t)^{-\frac{3}{2}(1-\frac{1}{p})},\, t\geq0,\,  2\leq p\leq\infty .\end{equation}
In 2016, Li and  Yong \cite{LY1} investigated the zero Mach number limit for the three-dimensional
model in the regime of smooth solutions. In  critical  Besov spaces, Danchin and Desjardins \cite{Dan6} and  Haspot \cite{HS3,HS4} obtained the  global well-posedness  of strong solutions close to a stable
equilibrium state. Charve et al. \cite{CDX} obtained the  global existence, Gevrey
analytic and algebraic time-decay estimates of   strong  solutions   when  the initial data are  close to a stable
equilibrium state in $L^{p}$-critical  framework. In 2019,   N. Chikami and T. Kobayashi \cite{NCTK} established the  global existence and algebraic time decay estimates of   strong  solutions   when  the initial data are  close to a stable
equilibrium state.

Due to the more complicated coupling between hyperbolic
equations and parabolic equations in the two-phase
flow model, the mathematical structure of the model is much more complex than that in the case of single phase
flow model. Therefore, extending  the currently available results for single phase flow
models  to two-phase models is not an easy task. Therefore,  more and more researchers  pay more attention to
 the mathematical problems of the generic two-phase model.  In \cite{BDG},  Bresch et al. first established the existence of global weak solutions to the 3D  generic two-fluid flow model \eqref{equ:CTFS}. Later Bresch-Huang-Li \cite{BHL} extended the result in \cite{BDG} and proved the existence of global weak solutions to the system \eqref{equ:CTFS} in one space dimension without capillarity terms.  In 2016, Evje-Wang-Wen \cite{EWW}   proved the global existence
of strong solutions to the  generic two-fluid flow  model \eqref{equ:CTFS} without capillary
terms    by the standard energy method under the condition that the initial data are
close to the constant equilibrium state in $H^{2}(\mathbb{R}^{3})$ and obtained the following optimal time decay rates for   strong solutions  if the initial data belong to $L^{1}$ additionally
\begin{equation}\label{eq:decate6}\|(\alpha^{+}\rho^{+}-1, u^{+}, \alpha^{-}\rho^{-}-1, u^{-})\|_{L^{2}}\leq C(1+t)^{-\frac{3}{4}},\, t\geq0,\end{equation}
\begin{equation}\label{eq:decate7}\|\nabla(\alpha^{+}\rho^{+}-1, u^{+}, \alpha^{-}\rho^{-}-1, u^{-})\|_{L^{2}}\leq C(1+t)^{-\frac{5}{4}},\, t\geq0.\end{equation}
Lai-Wen-Yao \cite{LWY}  studied the vanishing capillarity limit of  smooth solutions  to the  system  \eqref{equ:CTFS}   with  unequal pressure functions  if   $\|(\alpha^{+}\rho^{+}_0-1, \alpha^{-}\rho^{-}_0-1)\|_{H^{4}(\mathbb{R}^{3})} +\|(u^{+}_0, u^{-}_0)\|_{H^{3}(\mathbb{R}^{3})}$ are small enough.
Recently, based on complicated  spectral analysis  of  Green's function to the linearized system and on elaborate energy estimates
to the nonlinear system,  for the  system  \eqref{equ:CTFS},   authors \cite{CWYZ} proved  global solvability of  smooth solutions  close to an equilibrium state  in $H^s(\mathbb{R}^{3})(s\geq3)$ and further  got the  following time decay rates when the initial
perturbation is bounded in $L^{1}(\mathbb{R}^{3})$
\begin{equation}\label{eq:decate8}\|(\alpha^{+}\rho^{+}-1,  \alpha^{-}\rho^{-}-1)\|_{L^{2}}\leq C(1+t)^{-\frac{1}{4}},\, t\geq0,\end{equation}
\begin{equation}\label{eq:decate9}\|\nabla(\alpha^{+}\rho^{+}-1, \alpha^{-}\rho^{-}-1)\|_{L^{2}}\leq C(1+t)^{-\frac{3}{4}},\, t\geq0,\end{equation}
\begin{equation}\label{eq:decate10}\|( u^{+}, u^{-})\|_{L^{2}}\leq C(1+t)^{-\frac{3}{4}},\, t\geq0,\end{equation}
\begin{equation}\label{eq:decate11}\|\nabla( u^{+}, u^{-})\|_{L^{2}}\leq C(1+t)^{-\frac{5}{4}},\, t\geq0.\end{equation}
Here,  it should be pointed  out that  the functional spaces with high Sobolev regularity  is not \emph{the lowest index} in the sense of the scaling invariant of the associated System \eqref{equ:CTFS} and  the dimension of space is only  limited to $N=3$. Moreover, the decay rates of $(\alpha^{\pm}\rho^{\pm}-1)$ in \eqref{eq:decate8}-\eqref{eq:decate9}  are not the optimal.

The purpose of this work is to investigate the mathematical properties of the system \eqref{equ:CTFS} in critical regularity framework. More specifically, we address the question of whether available mathematical
results such as the global well-posedness and the optimal  time decay rates in critical Besov spaces to a single fluid governed by the compressible barotropic
Navier-Stokes equations may be extended to   multi-dimensional  non-conservative viscous
compressible two-fluid system. At this stage, we are going to use scaling considerations for the system \eqref{equ:CTFS} to
guess which spaces may be critical. One can check that if $(\alpha^{+}\rho^{+},\,u^{+},\,\alpha^{-}\rho^{-},\,u^{-})$  solves the system  \eqref{equ:CTFS}, so does $\big((\alpha^{+}\rho^{+})_{\lambda},\,u^{+}_{\lambda}, (\alpha^{-}\rho^{-})_{\lambda},\,u^{-}_{\lambda}\big)$ where:
\begin{equation}
\begin{split}\label{udivc}&(\alpha^{+}\rho^{+})_{\lambda}(t, x)=(\alpha^{+}\rho^{+})(\lambda^{2}t, \lambda x),\quad u^{+}_{\lambda}(t, x)=\lambda u^{+}(\lambda^{2}t, \lambda x),\\ &(\alpha^{-}\rho^{-})_{\lambda}(t, x)=(\alpha^{-}\rho^{-})(\lambda^{2}t, \lambda x),\quad u^{-}_{\lambda}(t, x)=\lambda u^{-}(\lambda^{2}t, \lambda x)
\end{split}
\end{equation}
provided that the pressure laws $P$ have been changed into $\lambda^{2}P$.  This suggests us to choose
initial data $\big((\alpha^{+}\rho^{+})_{0},\,u^{+}_{0},\,(\alpha^{-}\rho^{-})_{0},\,u^{-}_{0}\big)$ in critical  spaces whose norm is invariant for all $\lambda>0$ by the
transformation $\big((\alpha^{+}\rho^{+})_{0},\,u^{+}_{0},\,(\alpha^{-}\rho^{-})_{0},\,u^{-}_{0}\big)\rightarrow \big((\alpha^{+}\rho^{+})_{0}(\lambda\,\cdot),\,\lambda u^{+}_{0}(\lambda\,\cdot),\,(\alpha^{-}\rho^{-})_{0}(\lambda\,\cdot),\,\lambda u^{-}_{0}(\lambda\,\cdot)\big)$.  For the convenience of the reader, as in  \cite{BDG,CWYZ,LWY},  we also show some  derivations for
another expression of the pressure gradient in terms of the gradients of $\alpha^{+}\rho^{+}$ and $\alpha^{-}\rho^{-}$
 by using the pressure equilibrium assumption.
 The relation between the
pressures of the system \eqref{equ:CTFS} implies the following differential identities
\begin{equation}\label{1.3}\textrm{d}P^{+}=s_{+}^{2}\textrm{d}\rho^{+},\quad
\textrm{d}P^{-}=s_{-}^{2}\textrm{d}\rho^{-},\quad \hbox{where}\quad
s_{\pm}:=\sqrt{\frac{\textrm{d}P^{\pm}}{\textrm{d}\rho^{\pm}}(\rho^{\pm})}
=\sqrt{\overline{\gamma}^{\pm}\frac{P^{\pm}(\rho^{\pm})}{\rho^{\pm}}},
\end{equation}
where $s_{\pm}$ denote the sound speed of each phase respectively.

Let
\begin{equation}\label{1.4}
R^{\pm}=\alpha^{\pm}\rho^{\pm}.
\end{equation}
Resorting to $\eqref{equ:CTFS}_{1}$,  we have
\begin{equation}\label{1.5}
\textrm{d}\rho^{+}=\frac{1}{\alpha_{+}}(\textrm{d}R^{+}-
\rho^{+}\textrm{d}\alpha^{+}),
~~\textrm{d}\rho^{-}=\frac{1}{\alpha_{-}}(\textrm{d}R^{-}+
\rho^{-}\textrm{d}\alpha^{+}).
\end{equation}
Combining with \eqref{1.4} and \eqref{1.5}, we  conclude that
\begin{equation*}\label{1.6}
\textrm{d}\alpha^{+}=\frac{\alpha^{-}s_{+}^{2}}
{\alpha^{-}\rho^{+}s_{+}^{2}+\alpha^{+}\rho^{-}s_{-}^{2}}\textrm{d}R^{+}
-\frac{\alpha^{+}s_{-}^{2}}
{\alpha^{-}\rho^{+}s_{+}^{2}+\alpha^{+}\rho^{-}s_{-}^{2}}
\textrm{d}R^{-}.
\end{equation*}
Substituting  the above equality into \eqref{1.5}, we obtain
\begin{equation*}
\textrm{d}\rho^{+}=\frac{\rho^{+}\rho^{-}s_{-}^{2}}
{R^{-}(\rho^{+})^{2}s_{+}^{2}+R^{+}(\rho^{-})^{2}s_{-}^{2}}
\Big(\rho^{-}\textrm{d}R^{+}
+\rho^{+}\textrm{d}R^{-}\Big),
\end{equation*}
and
\begin{equation*}
\textrm{d}\rho^{-}=\frac{\rho^{+}\rho^{-}s_{+}^{2}}
{R^{-}(\rho^{+})^{2}s_{+}^{2}+R^{+}(\rho^{-})^{2}s_{-}^{2}}
\Big(\rho^{-}\textrm{d}R^{+}
+\rho^{+}\textrm{d}R^{-}\Big),
\end{equation*}
which give, for the pressure differential $\textrm{d}P^{\pm}$,
\begin{equation*}
\textrm{d}P^{+}=\mathcal{C}^{2}\big(\rho^{-}\textrm{d}R^{+}
+\rho^{+}\textrm{d}R^{-}\big),
\end{equation*}
and
\begin{equation*}
\textrm{d}P^{-}=\mathcal{C}^{2}\big(\rho^{-}\textrm{d}R^{+}
+\rho^{+}\textrm{d}R^{-}\big),
\end{equation*}
where
\begin{equation*}
\mathcal{C}^{2}\eqdefa\frac{s_{-}^{2}s_{+}^{2}}{\alpha^{-}\rho^{+}s_{+}^{2}
+\alpha^{+}\rho^{-}s_{-}^{2}}.
\end{equation*}
Recalling $\alpha^{+}+\alpha^{-}=1$, we get the following identity:
\begin{equation}\label{1.61}
\frac{R^{+}}{\rho^{+}}+\frac{R^{-}}{\rho^{-}}=1,\quad \hbox{and ~therefore~}
\rho^{-}=\frac{R^{-}\rho^{+}}{\rho^{+}-R^{+}}.
\end{equation}
Then it follows from the pressure relation $\eqref{equ:CTFS}_{4}$ that
\begin{equation}\label{1.7}
\varphi(\rho^{+}):=P^{+}(\rho^{+})-P^{-}(\frac{R^{-}\rho^{+}}
{\rho^{+}-R^{+}})=0.
\end{equation}
 Differentiating $\varphi$ with respect to $\rho^{+}$,   we have
$$\varphi^{'}(\rho^{+})=s^{2}_{+}+s^{2}_{-}\frac{R^{-}R^{+}}{(\rho^{+}-R^{+})^{2}}.$$
By the definition of $R^{+}$, it is natural to look for $\rho^{+}$ which belongs to $(R^{+},+\infty).$  Since $\varphi'>0$ in $(R^{+},+\infty)$  for any given $R^{\pm}>0,$ and $\varphi:(R^{+},+\infty)\longmapsto(-\infty,+\infty),$  this determines that $\rho^{+}=\rho^{+}(R^{+},R^{-})\in(R^{+},+\infty)$ is the unique solution of the equation \eqref{1.7}. Due to \eqref{1.5}, \eqref{1.61} and $\eqref{equ:CTFS}_{1}$, $\rho^{-}$ and $\alpha^{\pm}$  are defined as follows:
$$\rho^{-}(R^{+},R^{-})=\frac{R^{-}\rho^{+}(R^{+},R^{-})}{\rho^{+}(R^{+},R^{-})-R^{+}},$$
$$\alpha^{+}(R^{+},R^{-})=\frac{R^{+}}{\rho^{+}(R^{+},R^{-})},$$
$$\alpha^{-}(R^{+},R^{-})=1-\frac{R^{+}}{\rho^{+}(R^{+},R^{-})}=\frac{R^{-}}{\rho^{-}(R^{+},R^{-})}.$$
Based on  the above analysis,  the system \eqref{equ:CTFS} is equivalent to the following form
\begin{align}\label{equ:CTFS1}
\left\{
\begin{aligned}
&\p_tR^{\pm}+\textrm{div}(R^{\pm}u^{\pm})=0, \\
&\p_t(R^{+}u^{+})+\textrm{div}(R^{+}u^{+}\otimes u^{+})+\alpha^{+}\mathcal{C}^{2}[\rho^{-}\nabla R^{+}+\rho^{+}\nabla R^{-}]-\sigma^{+}R^{+}\nabla \Delta R^{+}
\\&\qquad=\textrm{div}\big(\alpha^{+}[\mu^{+}(\nabla u^{+}+\nabla^{t}u^{+})+\lambda^{+}\textrm{div}u^{+}\textrm{Id}]\big), \\
&\p_t(R^{-}u^{-})+\textrm{div}(R^{-}u^{-}\otimes u^{-})+\alpha^{-}\mathcal{C}^{2}[\rho^{-}\nabla R^{+}+\rho^{+}\nabla R^{-}]-\sigma^{-}R^{-}\nabla \Delta R^{-}
\\&\qquad=\textrm{div}\big(\alpha^{-}[\mu^{-}(\nabla u^{-}+\nabla^{t}u^{-})+\lambda^{-}\textrm{div}u^{-}\textrm{Id}]\big).
\end{aligned}
\right.
\end{align}
Here, we are concerned with the Cauchy problem of the system
 \eqref{equ:CTFS1} in $\mathbb{R}_{+}\times \mathbb{R}^N$ subject to the
initial data
\begin{equation}\label{eq:initial data}
(R^{+},\,u^{+},\,R^{-},\,u^{-})(x,t)|_{t=0}=(R^{+}_{0},\,u^{+}_{0},\,R^{-}_{0},\,u^{-}_{0})(x), \quad x\in\mathbb{R}^{N},
\end{equation}
and
$$u^{+}(x,t)\rightarrow0,\quad u^{-}(x,t)\rightarrow0,\quad R^{+}\rightarrow R^{+}_{\infty}>0,\quad R^{-}\rightarrow R^{-}_{\infty}>0, \hbox{~as~} |x|\rightarrow\infty,$$
where $R^{\pm}_{\infty}$ denote the background doping profile, and in the present paper $R^{\pm}_{\infty}$ are taken as $1$ without losing generality.

 For simplicity, we take $\sigma^{+}=\sigma^{-}=1$. Set
$c^{\pm}=R^{\pm}-1$.
Then, the system \eqref{equ:CTFS1} can be rewritten  as
\begin{equation}\label{equ:CTFS2}
\left\{
\begin{aligned}{}
&\p_tc^{+}+\textrm{div}u^{+}=H_{1},\\
&\p_t{u}^{+}+\beta_{1}\nabla c^{+}
+\beta_{2}\nabla c^{-}-\nu_{1}^{+}\Delta u^{+}
-\nu_{2}^{+}\nabla\textrm{div}u^{+}-\nabla \Delta c^{+}=H_{2},
\\&\p_tc^{-}+\textrm{div}u^{-}=H_{3},
\\&\p_tu^{-}+\beta_{3}\nabla c^{+}
+\beta_{4}\nabla c^{-}-\nu_{1}^{-}\Delta u^{-}
-\nu_{2}^{-}\nabla\textrm{div}u^{-}-\nabla \Delta c^{-}=H_{4},
\end{aligned}
\right.
\end{equation}
with initial data
\begin{equation}\label{equ:CTFS3}(c^{+},\,u^{+},\,c^{-},\,u^{-})|_{t=0}=(c^{+}_{0},\,u^{+}_{0},\,c^{-}_{0},\,u^{-}_{0}),\end{equation}
where $\beta_{1}=\frac{\mathcal{C}^{2}(1,1)\rho^{-}(1,1)}{\rho^{+}(1,1)},\quad
\beta_{2}=\beta_{3}=\mathcal{C}^{2}(1,1),\quad
\beta_{4}=\frac{\mathcal{C}^{2}(1,1)\rho^{+}(1,1)}{\rho^{-}(1,1)},\quad
\nu_{1}^{\pm}=\frac{\mu^{\pm}}{\rho^{\pm}(1,1)},\quad
\nu_{2}^{\pm}=\frac{\mu^{\pm}+\lambda^{\pm}}{\rho^{\pm}(1,1)}$
and the source terms are
\begin{align}\label{3.2}
H_{1}&=H_{1}(c^{+},u^{+})=-\textrm{div}(c^{+}u^{+}),\\
\label{3.3}
H_{2}^{i}&=H_{2}(c^{+},u^{+},c^{-})=-g_{+}(c^{+},c^{-})\partial_{i}c^{+}
-\tilde{g}_{+}(c^{+},c^{-})\partial_{i}c^{-}-(u^{+}\cdot\nabla)u_{i}^{+}\nonumber\\
&\quad+\mu^{+}h_{+}(c^{+},c^{-})\partial_{j}c^{+}\partial_{j}u^{+}_{i}
+\mu^{+}k_{+}(c^{+},c^{-})\partial_{j}c^{-}\partial_{j}u^{+}_{i}\nonumber\\
&\quad+\mu^{+}h_{+}(c^{+},c^{-})\partial_{j}c^{+}\partial_{i}u^{+}_{j}
+\mu^{+}k_{+}(c^{+},c^{-})\partial_{j}c^{-}\partial_{i}u^{+}_{j}\\
&\quad+\lambda^{+}h_{+}(c^{+},c^{-})\partial_{i}c^{+}\partial_{j}u^{+}_{j}
+\lambda^{+}k_{+}(c^{+},c^{-})\partial_{i}c^{-}\partial_{j}u^{+}_{j}\nonumber\\
&\quad+\mu^{+}l_{+}(c^{+},c^{-})\partial_{j}^{2}u_{i}^{+}+(\mu^{+}+\lambda^{+})l_{+}
(c^{+},c^{-})\partial_{i}\partial_{j}u^{+}_{j},\quad  i,j \in \{1,2,\cdots N\},\nonumber\\
\label{3.4}
H_{3}&=H_{3}(c^{-},u^{-})=-\textrm{div}(c^{-}u^{-}),\\
\label{3.5}
H_{4}^{i}&=H_{4}(c^{+},u^{-},c^{-})=-g_{-}(c^{+},c^{-})\partial_{i}c^{-}
-\tilde{g}_{-}(c^{+},c^{-})\partial_{i}c^{+}-(u^{-}\cdot\nabla)u_{i}^{-}\nonumber\\
&\quad+\mu^{-}h_{-}(c^{+},c^{-})\partial_{j}c^{+}\partial_{j}u^{-}_{i}
+\mu^{-}k_{-}(c^{+},c^{-})\partial_{j}c^{-}\partial_{j}u^{-}_{i}\nonumber\\
&\quad+\mu^{-}h_{-}(c^{+},c^{-})\partial_{j}c^{+}\partial_{i}u^{-}_{j}
+\mu^{-}k_{-}(c^{+},c^{-})\partial_{j}c^{-}\partial_{i}u^{-}_{j}\\
&\quad+\lambda^{-}h_{-}(c^{+},c^{-})\partial_{i}c^{+}\partial_{j}u^{-}_{j}
+\lambda^{-}k_{-}(c^{+},c^{-})\partial_{i}c^{-}\partial_{j}u^{-}_{j}\nonumber\\
&\quad+\mu^{-}l_{-}(c^{+},c^{-})\partial_{j}^{2}u_{i}^{-}+(\mu^{-}+\lambda^{-})l_{-}
(c^{+},c^{-})\partial_{i}\partial_{j}u^{-}_{j}, \quad  i,j \in \{1,2,\cdots N\},\nonumber
\end{align}
where we define the nonlinear functions of $(c^{+},c^{-})$ by
\begin{align}\label{3.6}
\left\{
\begin{aligned}
&g_{+}(c^{+},c^{-})=\frac{(\mathcal{C}^{2}\rho^{-})(c^{+}+1,c^{-}+1)}
{\rho^{+}(c^{+}+1,c^{-}+1)}
-\frac{(\mathcal{C}^{2}\rho^{-})(1,1)}
{\rho^{+}(1,1)},
\\&g_{-}(c^{+},c^{-})=\frac{(\mathcal{C}^{2}\rho^{+})(c^{+}+1,c^{-}+1)}
{\rho^{-}(c^{+}+1,c^{-}+1)}
-\frac{(\mathcal{C}^{2}\rho^{+})(1,1)}
{\rho^{-}(1,1)},
\end{aligned}
\right.
\end{align}
\begin{align}\label{3.7}
\left\{
\begin{aligned}
&h_{+}(c^{+},c^{-})=\frac{(\mathcal{C}^{2}\alpha^{-})(c^{+}+1,c^{-}+1)}
{[c^{+}+1]s_{-}^{2}(c^{+}+1,c^{-}+1)},
\\&h_{-}(c^{+},c^{-})=-\frac{\mathcal{C}^{2}(c^{+}+1,c^{-}+1)}
{(\rho^{-}s_{-}^{2})(c^{+}+1,c^{-}+1)},
\end{aligned}
\right.
\end{align}
\begin{align}\label{3.8}
\left\{
\begin{aligned}
&k_{+}(c^{+},c^{-})=-{\frac{\mathcal{C}^{2}(c^{+}+1,c^{-}+1)}
{[c^{-}+1](s_{+}^{2}\rho^{+})(c^{+}+1,c^{-}+1)}},
\\&k_{-}(c^{+},c^{-})={\frac{(\alpha^{+}\mathcal{C}^{2})(c^{+}+1,c^{-}+1)}
{[c^{-}+1]s_{+}^{2}(c^{+}+1,c^{-}+1)}},
\end{aligned}
\right.
\end{align}
\begin{align}\label{3.9}
\tilde{g}_{+}(c^{+},c^{-})=\tilde{g}_{-}(c^{+},c^{-})=\mathcal{C}^{2}(c^{+}+1,c^{-}+1)
-\mathcal{C}^{2}(1,1),
\end{align}
\begin{align}\label{3.10}
l_{\pm}(c^{+},c^{-})=\frac{1}{\rho^{\pm}(c^{+}+1,c^{-}+1)}-\frac{1}
{\rho^{\pm}(1,1)}.
\end{align}

Although the non-conservative viscous
compressible two-fluid  flow model, to some extent, is similar to the compressible Navier-Stokes equations,
it is non-trivial to apply directly the ideas used in single-phase models into the two-phase models since the momentum equation is given only for the mixture and that the pressure involves the masses of two phases in a nonlinear way. Now, let us explain on some of the main difficulties and techniques involved in the process. First,  the presences of the two-phase flows effect results in a different and more intricate coupling relations in the present system,  so that the overall analysis is considerably more
complicated. Second, in two-fluid  flow model,  we find that there have some  binary functions such as  $g_{\pm}(c^{+},c^{-})$ including  the two  variables $c^{+}$ and $c^{-}$ in \eqref{3.3} and \eqref{3.5}. In order to obtain the  estimates  of these nonlinear terms in Besov spaces, we need to exploit  continuity for the composition for binary functions (see  Lemma \ref{p27} for details), which are distinguished estimates for the complicated two-phase flow model. Third, employing the energy argument of Godunov
\cite{GO} for partially dissipative first-order symmetric systems (further developed by
 \cite{KF}), the Littlewood-Paley decomposition and  Fourier-Plancherel
theorem,  we can  obtain maximal regularity estimates for  the linearized  system in Besov spaces and show
the  parabolic properties of $(c^{+},\,u^{+},\,c^{-},\,u^{-})$ in all frequencies (see  Lemma \ref{all part} for details). Finally,  one may wonder  how global strong solutions constructed above look like for large time. Under a suitable additional condition involving only the low frequencies of the  data and in the  $L^{2}$critical regularity framework,  we  exhibit the optimal time decay rates  for  the  constructed global strong
solutions. In this part, our main ideas are based on  the low-high frequency decomposition   and  a refined time-weighted energy functional. In  low frequencies, making good use of
Fourier   localization  analysis to a linearized parabolic-hyperbolic system in order to obtain  smoothing effects of  Green's function  and  avoid some complicate   spectral analysis as in \cite{CWYZ}. In  high frequencies,  in order to close the energy estimates, we further exploit some  decay estimates with gain of
regularity   of $(c^{+},\,u^{+},\,c^{-},\,u^{-})$. With these  analysis tools in hand,  we finally establish the global well-posedness and the optimal time decay rates of strong solutions to the Cauchy problem \eqref{equ:CTFS1}-\eqref{eq:initial data}.

Now we state our main results as follows:
\begin{Theorem}\label{th:main1}  Assume that $\big((\sqrt{\beta_1}+\Lambda)(R^{+}_{0}-1),\,u^{+}_{0},\,(\sqrt{\beta_4}+\Lambda)(R^{-}_{0}-1),\,u^{-}_{0}\big)\in \dot{B}^{\f{N}{2}-1}_{2,1}\times
\dot{B}^{\f{N}{2}-1}_{2,1}\times \dot{B}^{\f{N}{2}-1}_{2,1}\times
\dot{B}^{\f{N}{2}-1}_{2,1}$. Then there exists a constant $\eta>0$ such that if
\begin{equation}
X(0)\eqdefa\big\|\big((\sqrt{\beta_1}+\Lambda)(R^{+}_{0}-1),u^{+}_0,(\sqrt{\beta_4}+\Lambda)(R^{-}_{0}-1),u^{-}_0\big)\big\|_{
\dot{B}^{\f{N}{2}-1}_{2,1}}\leq \eta,  \label{1Hn/2}
\end{equation}%
then the Cauchy  problem \eqref{equ:CTFS1}-\eqref{eq:initial data} admits a unique global solution $(R^{+}-1,\,u^{+},\,R^{-}-1,\,u^{-})$ satisfying that for all $t\geq 0$,
\begin{equation}\label{1.6A}
\begin{split}
\,X(t)\lesssim X(0),
\end{split}
\end{equation}
where
\begin{equation}
\begin{split}
\label{1.6}
\,X(t)&\eqdefa\big\|\big((\sqrt{\beta_1}+\Lambda)(R^{+}-1),u^{+},(\sqrt{\beta_4}+\Lambda)(R^{-}-1),u^{-}\big)\big\|_{\wt L^\infty_t(\dot B^{\frac{N}{2}-1}_{2,1})}
\\&\quad+\big\|\big((\sqrt{\beta_1}+\Lambda)(R^{+}-1),u^{+},(\sqrt{\beta_4}+\Lambda)(R^{-}-1),u^{-}\big)\big\|_{L^1_t(\dot B^{\frac{N}{2}+1}_{2,1})}.
\end{split}
\end{equation}
\end{Theorem}
\begin{Theorem}\label{th:decay}
Let the data $(R^{+}_{0}-1,\,u^{+}_{0},\,R^{-}_{0}-1,\,u^{-}_{0})$ satisfy the assumptions of Theorem \ref{th:main1}. Denote $\langle \tau\rangle\eqdefa\sqrt{1+\tau^2}$
and $\alpha\eqdefa\frac{N}{2}+\frac{1}{2}-\varepsilon$ with $\varepsilon>0$ arbitrarily small.
 There exists a positive constant $c$ such that if in addition
\begin{equation}\label{eq:D0}
D_0\eqdefa\big\|\big((\sqrt{\beta_1}+\Lambda)(R^{+}_{0}-1),u^{+}_0,(\sqrt{\beta_4}+\Lambda)(R^{-}_{0}-1),u^{-}_0\big)\big\|^{\ell}_{\dot B^{-\frac N2}_{2,\infty}}\leq c,
\end{equation}
then  the global solution $(R^{+}-1,\,u^{+},\,R^{-}-1,\,u^{-})$
given by Theorem \ref{th:main1} satisfies for all $t\geq0,$
\begin{equation}
\label{1.8}
D(t)\leq C\Big(D_0+\big\|\big(\Lambda R^{+}_{0},\,u^{+}_{0},\,\Lambda R^{-}_{0},\,u^{-}_{0}\big)\big\|^h_{\dot B^{\frac N2-1}_{2,1}}\Big)
\end{equation}
with
\begin{equation}
\begin{split}
\label{1.9}
D(t)&\eqdefa \sup_{s\in(\varepsilon-\frac N2,2]}\big\|\langle\tau\rangle^{\frac N4+\frac s2}\big((\sqrt{\beta_1}+\Lambda)(R^{+}-1),u^{+},(\sqrt{\beta_4}+\Lambda)(R^{-}-1),u^{-}\big)\big\|_{L^\infty_t(\dot B^s_{2,1})}^\ell
\\ &\quad+\big\|\tau^{\alpha}\Lambda^{2}\big((\sqrt{\beta_1}+\Lambda)(R^{+}-1),u^{+},(\sqrt{\beta_4}+\Lambda)(R^{-}-1),u^{-}\big)\big\|_{\wt L^\infty_t(\dot B^{\frac N2-1}_{2,1})}^h.
\end{split}
\end{equation}
\end{Theorem}
\begin{Remark}\label{1.2} Compared with \cite{CWYZ,EWW},   in theorems \ref{th:main1} and \ref{th:decay}, we obtain the global well-posedness and the     optimal time decay rates for multi-dimensional  non-conservative viscous
compressible two-fluid system
 \eqref{equ:CTFS}  in critical regularity framework respectively. Additionally, in Theorems \ref{th:decay},
the regularity  index $s$ can take both negative and nonnegative values, rather than only nonnegative integers, which   improves the classical decay results  in high Sobolev regularity, such as \cite{CWYZ,EWW,Tan1,Tan2}.
 \end{Remark}
\begin{Remark}\label{1.2}Due to the embedding $L^{1}(\mathbb{R}^3)\hookrightarrow \dot{B}^{-\frac 32}_{2,\infty}(\mathbb{R}^3)$,  our results in Theorem \ref{th:decay} extend the known conclusions in \cite{CWYZ,EWW,Tan1,Tan2}. In particular, our condition involves only the low frequencies of the data and is based on the  $L^{2}(\mathbb{R}^{N})$-norm framework.  In particular, the decay rates of strong solutions is
in the so-called critical Besov spaces in any dimension $N \geq 2$ and the dimension of space is more extensive and is not limited to $N=3$.
\end{Remark}
\begin{Remark} It should be pointed out that we only deal with  the case for $L^{2}(\mathbb{R}^{N})$-framework in this paper. The  case for $L^{p}(\mathbb{R}^{N})$-framework will be considered in our future work.
\end{Remark}

The rest of the paper unfolds as follows. In the next section,  we recall some basic facts about Littlewood-Paley
decomposition, Besov spaces  and some useful lemmas.
 In Section 3, we will exploit
 maximal regularity estimates for  the linearized  system in Besov Spaces.
 Section 4 is devoted to the proof of the global well-posedness for initial data near equilibrium in critical Besov spaces.
In Section 5,  we present  the optimal time decay rates of the  global  strong solutions. At last, in Section 6,   we show  some corollaries  and the standard optimal   $L^{q}$-$L^{r}$ time decay rates.

\noindent{\bf Notations.}  We assume $C$ be
a positive generic constant throughout this paper that may vary at
different places and
denote  $A\le CB$ by  $A\lesssim B$.
We shall also use the  following notations
 $$z^\ell\eqdefa\sum_{j\leq k_0}\dot{\Delta}_{j}z\quad\hbox{and}\quad z^h\eqdefa z-z^\ell, \quad\hbox{for some } j_0.$$
$$\|z\|^\ell_{\dot B^s_{2,1}}\eqdefa \sum_{j\leq k_0}2^{js}\|\dot{\Delta}_{j}z\|_{L^{2}}\quad\hbox{and}\quad \|z\|^h_{\dot B^s_{2,1}}\eqdefa \sum_{j\geq k_0}2^{js}\|\dot{\Delta}_{j}z\|_{L^{2}}, \quad\hbox{for some } j_0.$$
Noting the small overlap between low and high frequencies,  we have
$$\|z^\ell\|_{\dot B^s_{2,1}}\lesssim \|z\|^\ell_{\dot B^s_{2,1}}\quad\hbox{and}\quad \|z^h\|_{\dot B^s_{2,1}}\lesssim \|z\|^h_{\dot B^s_{2,1}}.$$
\par

\section{ Littlewood-Paley Theory and Some Useful Lemmas }
\ \ \ \ \ Let us introduce the Littlewood-Paley decomposition.
Choose a radial function  $\varphi \in {\cS}(\mathbb{R}^N)$
supported in ${\cC}=\{\xi\in\mathbb{R}^N,\,
\frac{3}{4}\le|\xi|\le\frac{8}{3}\}$ such that \beno \sum_{q\in
\mathbb{Z}}\varphi(2^{-q}\xi)=1 \quad \textrm{for all}\,\,\xi\neq 0.
\eeno The homogeneous frequency localization operators
$\dot{\Delta}_q$ and $\dot{S}_q$ are defined by
\begin{align}
\dot{\Delta}_qf=\varphi(2^{-q}D)f,\quad \dot{S}_qf=\sum_{k\le
q-1}\dot{\Delta}_kf\quad\mbox{for}\quad q\in \mathbb{Z}. \nonumber
\end{align}
With our choice of $\varphi$, one can easily verify that
\begin{equation*}\begin{split}
&\dot{\Delta}_q\dot{\Delta}_kf=0\quad \textrm{if}\quad|q-k|\ge
2\quad \textrm{and}
\quad \\
&\dot{\Delta}_q(\dot{S}_{k-1}f\dot{\Delta}_k f)=0\quad
\textrm{if}\quad|q-k|\ge 5.
\end{split}\end{equation*}
 We denote the space ${\cZ'}(\mathbb{R}^N)$ by the dual space of
${\cZ}(\mathbb{R}^N)=\{f\in {\cS}(\mathbb{R}^N);\,D^\alpha
\hat{f}(0)=0; \forall\alpha\in\mathbb{ N}^N \,\mbox{multi-index}\}.$
It also can be identified by the quotient space of
${\cS'}(\mathbb{R}^N)/{\cP}$ with the polynomials space ${\cP}$. The
formal equality \beno f=\sum_{q\in\mathbb{Z}}\dot{\Delta}_qf \eeno
holds true for $f\in {\cZ'}(\mathbb{R}^N)$ and is called the
homogeneous Littlewood-Paley decomposition.

The following Bernstein's inequalities will be frequently used.
\begin{Lemma}\cite{Che-book}\label{Lem:Bernstein}
Let $1\le p_{1}\le p_{2}\le+\infty$. Assume that $f\in L^{p_{1}}(\mathbb{R}^N)$,
then for any $\gamma\in(\mathbb{N}\cup\{0\})^N$, there exist
constants $C_1$, $C_2$ independent of $f$, $q$ such that \beno
&&{\rm supp}\hat f\subseteq \{|\xi|\le A_02^{q}\}\Rightarrow
\|\partial^\gamma f\|_{p_2}\le C_12^{q{|\gamma|}+q
N(\frac{1}{p_1}-\frac{1}{p_2})}\|f\|_{p_1},
\\
&&{\rm supp}\hat f\subseteq \{A_12^{q}\le|\xi|\le
A_22^{q}\}\Rightarrow \|f\|_{p_1}\le
C_22^{-q|\gamma|}\sup_{|\beta|=|\gamma|}\|\partial^\beta f\|_{p_1}.
\eeno
\end{Lemma}
Let us recall the definition of homogeneous
Besov spaces  (see \cite{BCD,Dan1}).
\begin{Definition}\label{def2.2} Let $s\in \mathbb{R}$, $1\le p,
r\le+\infty$. The homogeneous Besov space $\dot{B}^{s}_{p,r}$ is
defined by
$$\dot{B}^{s}_{p,r}=\Big\{f\in {\cZ'}(\mathbb{R}^N):\,\|f\|_{\dot{B}^{s}_{p,r}}<+\infty\Big\},$$
where \beno \|f\|_{\dot{B}^{s}_{p,r}}\eqdefa \Big\|2^{qs}
\|\dot{\Delta}_qf(t)\|_{p}\Big\|_{\ell^r}.\eeno
\end{Definition}
\begin{Remark}\label{2.3}\quad
Some properties about the  Besov spaces are as follows
\begin{itemize}
\item\,\, Derivation: $$\|f\|_{\dot{B}^{s}_{2,1}}\approx\|\nabla f\|_{\dot{B}^{s-1}_{2,1}};$$
\item\,\, Algebraic properties: for $s > 0$, $\dot{B}^{s}_{2,1} \cap L^{\infty}
$ is an algebra;
\item\,\,Interpolation: for
$s_1, s_2\in\mathbb{R}$ and $\theta\in[0,1]$,
we have $$\|f\|_{\dot{B}^{\theta s_1+(1-\theta)s_2}_{2,1}}\le
\|f\|^\theta_{\dot{B}^{s_1}_{2,1}}\|f\|^{(1-\theta)}_{\dot{B}^{s_2}_{2,1}}.$$
\end{itemize}
\end{Remark}
\begin{Definition} Let $s\in \mathbb{R}$, $1\le p,\rho,
r\le+\infty$. The homogeneous space-time  Besov space
$L^\rho_{T}(\dot{B}_{p,  r}^s)$ is defined by
$$L^\rho_{T}(\dot{B}_{p,  r}^s)=\Big\{f\in \mathbb{R}_{+}\times{\cZ'}(\mathbb{R}^N):\,\|f\|_{L^\rho_{T}(\dot{B}_{p,  r}^s)} <+\infty\Big\},$$ where $$
\|f\|_{L^\rho_{T}(\dot{B}_{p,  r}^s)}\eqdefa \Big\|\big\| 2^{qs}
\|\dot{\Delta}_q f\|_{L^p}\big\|_{\ell ^{r}}\Big\|_{L^\rho_{T}}.$$
\end{Definition}
We next introduce the Besov-Chemin-Lerner space
$\widetilde{L}^q_T(\dot{B}^{s}_{p,r})$ which is initiated in
\cite{Che-Ler}.
\begin{Definition}Let $s\in \mathbb{R}$, $1\le
p,q,r\le+\infty$, $0<T\le+\infty$. The space
$\widetilde{L}^q_T(\dot{B}^s_{p,r})$ is defined by
$$\widetilde{L}^q_T(\dot{B}^s_{p,r})=\Big\{f\in \mathbb{R}_{+}\times{\cZ'}(\mathbb{R}^N):\,\|f\|_{\widetilde{L}^q_T(\dot{B}^{s}_{p,r})}<+\infty\Big\},$$
where
$$\|f\|_{\widetilde{L}^q_T(\dot{B}^{s}_{p,r})}\eqdefa \Bigl\|2^{qs}
\|\dot{\Delta}_qf(t)\|_{L^q(0,T;L^p)}\Bigr\|_{\ell^r}.$$
\end{Definition}
Obviously, $
\widetilde{L}^1_T(\dot{B}^s_{p,1})=L^1_T(\dot{B}^s_{p,1}). $ By  a
direct application  of  Minkowski's inequality, we have the
following relations between these spaces
\begin{equation*}
L^\rho_{T}(\dot{B}_{p,r}^s)\hookrightarrow\widetilde
L^\rho_{T}(\dot{B}_{p,r}^s),\,\textnormal{if}\quad  r\geq
\rho,\end{equation*}
\begin{equation*}
\widetilde L^\rho_{T}(\dot{B}_{p,r}^s)\hookrightarrow
L^\rho_{T}(\dot{B}_{p,r}^s),\, \textnormal{if}\quad \rho\geq r.
\end{equation*}
The usual product is continuous in many Besov spaces. The following proposition
the proof of which may be found in \cite{RS} section 4.4,
see in particular inequality (28) page 174)  will be very useful.
\begin{Proposition}\label{p26}
 For all $1\leq r,p, p_1, p_2\leq+\infty$,  there exists a positive universal
 constant such that
$$\|fg\|_{\dot{B}^{s}_{p,r}}\lesssim
\|f\|_{L^\infty}\|g\|_{\dot{B}^{s}_{p,r}}+\|g\|_{L^\infty}\|f\|_{\dot{B}^{s}_{p,r}},
\quad \text{if}\quad s>0;$$
$$\|fg\|_{\dot{B}^{s_1+s_2-\frac{N}{p}}_{p,r}}\lesssim
\|f\|_{\dot{B}^{s_1}_{p,r}}\|g\|_{\dot{B}^{s_2}_{p,\infty}}, \quad
\text{if}\quad s_1,s_2<\frac{N}{p},\quad \text{and}\quad
s_1+s_2>0;$$
$$\|fg\|_{\dot{B}^{s}_{p,r}}\lesssim
\|f\|_{\dot{B}^{s}_{p,r}}\|g\|_{\dot{B}^{\frac{N}{p}}_{p,\infty}\cap
L^{\infty}}, \quad \text{if}\quad |s|<\frac{N}{p};$$
$$\|fg\|_{\dot{B}^s_{2,1}}\lesssim \|f\|_{\dot{B}^{N/2}_{2,1}}\|g\|_{\dot{B}^s_{2,1}}, \quad \text{if}\quad
s\in (-N/2,N/2].$$
\end{Proposition}
For the composition of the  binary functions, we have the following estimates.
\begin{Lemma}\label{p27}
\mbox{(i) } Let $s>0$, $t\geq0$, $1\le p, q,r\le \infty$ and $(f,g)\in \big(\tilde{L}_{t}^{q}(\dot{B}^{s}_{p,r})\cap L_{t}^{\infty}(L^{\infty})\big)\times\big(\tilde{L}_{t}^{q}(\dot{B}^{s}_{p,r})\cap L_{t}^{\infty}(L^{\infty})\big)$. If $F\in W_{loc}^{[s]+2,\infty}(\mathbb{R})\times W_{loc}^{[s]+2,\infty}(\mathbb{R})$ with
$F(0,0)=0$, then $F(f,g)\in \tilde{L}_{t}^{q}(\dot{B}^{s}_{p,r})$. Moreover,  there exists a
 $C$ depending only on $s,p,N$ and $F$, and
such that
\begin{equation}\label{eq:2.1}
\|F(f,g)\|_{\tilde{L}_{t}^{q}(\dot{B}^{s}_{p,r})}
\leq C\big(1+\|f\|_{L_{t}^{\infty}(L^{\infty})}\|g\|_{L_{t}^{\infty}(L^{\infty})}\big)^{[s]+1}
\|(f,g)\|_{\tilde{L}_{t}^{q}(\dot{B}^{s}_{p,r})}.
\end{equation}
\mbox{ (ii)} If $(f_1,g_1)\in \tilde{L}_{t}^{\infty}(\dot{B}^{\frac{N}{p}}_{p,1})\times\tilde{L}_{t}^{\infty}(\dot{B}^{\frac{N}{p}}_{p,1})$ and $(f_2,g_2)\in \tilde{L}_{t}^{\infty}(\dot{B}^{\frac{N}{p}}_{p,1})\times\tilde{L}_{t}^{\infty}(\dot{B}^{\frac{N}{p}}_{p,1})$, and  $(f_2-f_1,g_2-g_1)$ belongs to $\tilde{L}_{t}^{q}(\dot{B}^{s}_{p,r})\times \tilde{L}_{t}^{q}(\dot{B}^{s}_{p,r})$ with $s\in (-\frac{N}{p},\frac{N}{p}]$ . If $F\in W_{loc}^{[\frac{N}{p}]+2,\infty}(\mathbb{R})\times W_{loc}^{[\frac{N}{p}]+2,\infty}(\mathbb{R})$ with
$\partial_{1}F(0,0)=0$ and $\partial_{2}F(0,0)=0$. Then, there exists a
 $C$ depending only on $s, p, N$ and $F$
such that
\begin{equation}\label{eq:2.3}\begin{split}
&\|F(f_2,g_2)-F(f_1,g_1)\|_{\tilde{L}_{t}^{q}(\dot{B}^{s}_{p,r})}
\\&\leq C\Big(1+\|(f_1,f_2)\|_{L_{t}^{\infty}(L^{\infty})}\|(g_1,g_2)\|_{L_{t}^{\infty}(L^{\infty})}\Big)^{[\frac{N}{p}]+1} \Big(\|(f_1,g_1)\|_{\tilde{L}_{t}^{\infty}(\dot{B}^{\frac{N}{p}}_{p,1})} +\|(f_2,g_2)\|_{\tilde{L}_{t}^{\infty}(\dot{B}^{\frac{N}{p}}_{p,1})}\Big)
\\&\quad\times\Big(\|f_2-f_1\|_{\tilde{L}_{t}^{q}(\dot{B}^{s}_{p,r})}+\|g_2-g_1\|_{\tilde{L}_{t}^{q}(\dot{B}^{s}_{p,r})}\Big).
\end{split}
\end{equation}
\end{Lemma}

\noindent{\bf Proof.} \mbox{(i)}
From  the smoothness of  $F$  and $F(0,0)=0,$  we have  the following formal decomposition
\begin{equation*}
\begin{split}
F(f,g)&=\sum_{j\in\mathbb{Z}}
[F(\dot{S}_{j+1}f,\dot{S}_{j+1}g)-F(\dot{S}_{j}f,\dot{S}_{j}g)]
\\&=\sum_{j\in\mathbb{Z}}
 \int_{0}^{1}\partial_{1}F(\dot{S}_{j}f+\theta \dot{\Delta}_{j}f,\dot{S}_{j}g+\theta\dot{\Delta}_{j}g)d\theta \dot{\Delta}_{j}f
\\&\quad+\int_{0}^{1}\partial_{2} F(\dot{S}_{j}f+\theta \dot{\Delta}_{j}f,\dot{S}_{j}g+\theta\dot{\Delta}_{j}g)d\theta \dot{\Delta}_{j}g\\
&\triangleq\sum\limits_{j}m_{j}^{1}\dot{\Delta}_{j}f+m_{j}^{2}\dot{\Delta}_{j}g.
\end{split}
\end{equation*}
First, we have the following claim:
\begin{equation}\label{equ:1.2}
\|\partial^{\alpha} m_{j}^{i}\|_{L^{\infty}}\leq C\big(1+\|f\|_{L^{\infty}}\|g\|_{L^{\infty}}\big)^{|\alpha|}2^{j|\alpha|}\|\partial_iF\|_{W^{|\alpha|,\infty}},\quad \forall\alpha\in\mathbb{N}^{N}, ~i=1,2.
\end{equation}
Indeed, taking $m=2, h^{1}=\dot{S}_{j}f+\theta \dot{\Delta}_{j}f, h^{2}=\dot{S}_{j}g+\theta\dot{\Delta}_{j}g,$
$h_{\alpha}=(h_{\alpha}^{1},h_{\alpha}^{2})
=(\partial_{x}^{\alpha}h^{1},\partial_{x}^{\alpha}h^{2}),$  from multivariate Fa\'{a} di Bruno formula \cite{CS}, we have
\begin{equation*}
\begin{split}
\partial^{\alpha}m_{j}^{i}
&=\int_{0}^{1}\partial^{\alpha} \partial_{i} F(\dot{S}_{j}f+\theta \partial_{j}f,\dot{S}_{j}g+\theta\dot{\Delta}_{j}g)d\theta
=\int_{0}^{1}\partial^{\alpha} \partial_{i} F(h^{1},h^{2})d\theta
\\&=\int_{0}^{1} \sum\limits_{1\leq|\lambda|\leq |\alpha|}\partial^{\lambda}_{x}\partial_{i} F\sum_{s=1}^{n}
\sum\limits_{p_{s}(\alpha,\lambda)}(\alpha!)
\prod\limits_{q=1}^{s}
\frac{[h_{\ell_{q}}]^{k_{q}}}{(k_{q}!)[\ell_{q}!]^{|k_{q}|}}d\theta,
\end{split}
\end{equation*}
where $k_{q}=(k_{q_{1}},k_{q_{2}}),
\sum_{q=1}^{s}k_{q}=\lambda,
\sum_{q=1}^{s}|k_{q}|l_{q}=\alpha$,  $[h_{\ell_{q}}]^{k_{q}}=[\partial^{\ell_{q}}h^{1}]^{k_{q_{1}}}[\partial^{\ell_{q}}h^{2}]^{k_{q_{2}}}$.
Using H\"{o}lder's inequality and Bernstein's inequality,  we have
\begin{equation*}
\begin{split}
\|[h_{\ell_{q}}]^{k_{q}}\|_{L^{\infty}}
&\leq\|[\partial^{\ell_{q}}h^{1}]^{k_{q_{1}}}
\|_{L^{\infty}}
\|[\partial^{\ell_{q}}h^{2}]^{k_{q_{2}}}
\|_{L^{\infty}}
\\&\leq\|\partial^{\ell_{q}}(\dot{S}_{j}f+\theta \dot{\Delta}_{j}f)
\|_{L^{\infty}}^{k_{q_{1}}}
\|\partial^{\ell_{q}}(\dot{S}_{j}g+\theta\dot{\Delta}_{j}g)
\|_{L^{\infty}}^{k_{q_{2}}}
\\&\leq2^{j|\ell_{q}|k_{q_{1}}}\|\dot{S}_{j}f+\theta \dot{\Delta}_{j}f
\|_{L^{\infty}}^{k_{q_{1}}}
2^{j|\ell_{q}|k_{q_{2}}}\|\dot{S}_{j}g+\theta\dot{\Delta}_{j}g
\|_{L^{\infty}}^{k_{q_{2}}}
\\&\leq2^{j|\ell_{q}|(k_{q_{1}}+k_{q_{2}})}\|\dot{S}_{j}f+\theta \dot{\Delta}_{j}f
\|_{L^{\infty}}^{k_{q_{1}}}
\|(\dot{S}_{j}g+\theta\dot{\Delta}_{j}g)
\|^{k_{q_{2}}}_{L^{\infty}}
\\&\leq C2^{j|\ell_{q}||k_{q}|}\|f\|_{L^{\infty}}^{k_{q_1}}
\|g\|_{L^{\infty}}^{k_{q_2}}.
\end{split}
\end{equation*}
Thus,
\begin{equation*}
\begin{split}
\|\partial^{\alpha}m_{j}^{i}\|_{L^{\infty}}
&\leq C\prod\limits_{q=1}^{|\alpha|}2^{j|\ell_{q}||k_{q}|}
\|u\|_{L^{\infty}}^{k_{q_1}}\|v\|_{L^{\infty}}^{k_{q_2}}\|\partial_iF\|_{W^{|\alpha|,\infty}}
\\&\leq C\big(1+\|f\|_{L^{\infty}}\|g\|_{L^{\infty}}\big)^{|\alpha|}2^{j|\alpha|}\|\partial_iF\|_{W^{|\alpha|,\infty}},
\end{split}
\end{equation*}
where we have used
\begin{equation*}
\begin{split}
\sum_{q=1}^{|\alpha|}|\ell_{q}||k_{q}|=|\alpha|.
\end{split}
\end{equation*}
Next, we prove that $F(f,g)$ belongs to $\tilde{L}^{q}(\dot{B}^{s}_{p,r}).$  For all $j'\in \mathbb{Z}$, we have
\begin{equation*}
\begin{split}
\dot{\Delta}_{j'}F&=\sum_{j\leq j'}\dot{\Delta}_{j'}(m_{j}^{1}\dot{\Delta}_{j}f)
+\sum_{j\geq j'}\dot{\Delta}_{j'}(m_{j}^{1}\dot{\Delta}_{j}f)
+\sum_{j\leq j'}\dot{\Delta}_{j'}(m_{j}^{2}\dot{\Delta}_{j}g)
+\sum_{j\geq j'}\dot{\Delta}_{j'}(m_{j}^{2}\dot{\Delta}_{j}g)
\\&\triangleq\dot{\Delta}1_{j'}^{l}+\dot{\Delta}1_{j'}^{h}+\dot{\Delta}2_{j'}^{l}+\dot{\Delta}2_{j'}^{h}.
\end{split}
\end{equation*}
For the first term  $\dot{\Delta}1_{j'}^{l},$  taking advantage of Bernstein's inequality,  Leibniz's formula and \eqref{equ:1.2},  we have
\begin{equation*}
\begin{split}
2^{j's}\|\dot{\Delta}1_{j'}^{l}\|_{L_{t}^{q}(L^{p})}&\leq\sum_{j\leq j'}2^{j's}\|\dot{\Delta}_{j'}(m_{j}^{1}\Delta_{j}f)\|_{L_{t}^{q}(L^{p})}\\
&\leq C\sum_{j\leq j'}2^{(s-|\sigma|)j'}\|\partial^{\sigma} \dot{\Delta}_{j'}(m_{j}^{1}\dot{\Delta}_{j}f)\|_{L_{t}^{q}(L^{p})}
\\&\leq C\sum_{j\leq j'}2^{(s-|\sigma|)j'}C^{\beta}_{\sigma}\|\partial^{\beta}m_{j}^{1}\|_{L_{t}^{\infty}(L^{\infty})}
\|\partial^{\sigma-\beta}\dot{\Delta}_{j}f\|_{L_{t}^{q}(L^{p})}
\\&\leq C\sum_{j\leq j'}2^{(s-|\sigma|)j'}
\big(1+\|f\|_{L_{t}^{\infty}(L^{\infty})}\|g\|_{L_{t}^{\infty}(L^{\infty})}\big)^{|\beta|}2^{j|\beta|}\|\partial_1F\|_{W^{|\beta|,\infty}}
2^{j(|\sigma|-|\beta|)}\|\dot{\Delta}_{j}f\|_{L_{t}^{q}(L^{p})}
\\&\leq C\|\partial_1F\|_{W^{|\sigma|,\infty}}
\big(1+\|f\|_{L_{t}^{\infty}(L^{\infty})}\|g\|_{L_{t}^{\infty}(L^{\infty})}\big)^{|\sigma|}\sum_{j\leq j'}
2^{(s-|\sigma|)(j'-j)}\|\dot{\Delta}_{j}f\|_{L_{t}^{q}(L^{p})}2^{js}.
\end{split}
\end{equation*}
Taking $|\sigma|=[s]+1,$ convolution inequalities enable us to conclude that
$$(\sum_{j'}2^{j'sr}\|\dot{\Delta}1_{j'}^{l}\|^{r}_{L_{t}^{q}(L^{p})})^{\frac{1}{r}}
\leq  C\|\partial_1F\|_{W^{|\sigma|,\infty}}
\big(1+\|f\|_{L_{t}^{\infty}(L^{\infty})}\|g\|_{L_{t}^{\infty}(L^{\infty})}\big)^{[s]+1}
\|f\|_{\tilde{L}_{t}^{q}(\dot{B}^{s}_{p,r})}.$$
Similarly,
$$(\sum_{j'}2^{j'sr}\|\dot{\Delta}2_{j'}^{l}\|^{r}_{L_{t}^{q}(L^{p})})^{\frac{1}{r}}
\leq  C\|\partial_2F\|_{W^{|\sigma|,\infty}}
\big(1+\|f\|_{L_{t}^{\infty}(L^{\infty})}\|g\|_{L_{t}^{\infty}(L^{\infty})}\big)^{[s]+1}
\|g\|_{\tilde{L}_{t}^{q}(\dot{B}^{s}_{p,r})}.$$
Bounding the term pertaining to $\dot{\Delta}1_{j'}^{h}$
is easy. Indeed, we have according to  H\"{o}lder's inequality \eqref{equ:1.2},
\begin{equation*}
\begin{split}
2^{j's}\|\dot{\Delta}1_{j'}^{h}\|_{L_{t}^{q}(L^{p})}&\leq \sum_{j'\leq j}2^{j's}\|\dot{\Delta}_{j'}(m_{j}^{1}\dot{\Delta}_{j}f)\|_{L_{t}^{q}(L^{p})}\\
&\leq C \sum_{j'\leq j}2^{j's}\|m_{j}^{1}\|_{L_{t}^{\infty}(L^{\infty})}\|\dot{\Delta}_{j}f\|_{L_{t}^{q}(L^{p})}
\\&\leq C
\|\partial_1F\|_{L^{\infty}}\sum_{j'\leq j}2^{(j'-j)s}2^{js}\|\dot{\Delta}_{j}f\|_{L_{t}^{q}(L^{p})}.
\end{split}
\end{equation*}
Due to $s>0,$ we have
$$(\sum_{j'}2^{j'sr}\|\dot{\Delta}1_{j'}^{h}\|^{r}_{L_{t}^{q}(L^{p})})^{\frac{1}{r}}
\leq C
\|\partial_1F\|_{L^{\infty}}\|f\|_{\tilde{L}_{t}^{q}(\dot{B}^{s}_{p,r})}.$$
Similarity,
$$(\sum_{j'}2^{j'sr}\|\dot{\Delta}2_{j'}^{h}\|^{r}_{L_{t}^{q}(L^{p})})^{\frac{1}{r}}
\leq C
\|\partial_2F\|_{L^{\infty}}\|g\|_{\tilde{L}_{t}^{q}(\dot{B}^{s}_{p,r})}.$$
Thus, we deduce  \eqref{eq:2.1}.

\mbox{ (ii)}  We use the following identity
\begin{equation*}
\begin{split}
F(f_{2},g_2)-F(f_{1},g_1)&=
 \int_{0}^{1}\partial_{1}F(f_{1}+\theta(f_{2}-f_1),g_{1}+\theta(g_{2}-g_1))d\theta (f_{2}-f_1)
\\&\quad+\int_{0}^{1}\partial_{2} F(f_{1}+\theta(f_{2}-f_1),g_{1}+\theta(g_{2}-g_1))d\theta (g_{2}-g_1)\\
&\triangleq M_1(f_{2}-f_1)+M_2(g_{2}-g_1).
\end{split}
\end{equation*}
From Lemma \ref{p26} we obtain
\begin{equation*}
\begin{split}
\|F(f_{2},g_2)-F(f_{1},g_1)\|_{\tilde{L}_{t}^{q}(\dot{B}^{s}_{p,r})}&\leq \|M_1\|_{\tilde{L}_{t}^{\infty}(\dot{B}^{\frac{N}{p}}_{p,1})} \|f_{2}-f_1\|_{\tilde{L}_{t}^{q}(\dot{B}^{s}_{p,r})}+\|M_2\|_{\tilde{L}_{t}^{\infty}(\dot{B}^{\frac{N}{p}}_{p,1})} \|g_{2}-g_1\|_{\tilde{L}_{t}^{q}(\dot{B}^{s}_{p,r})}.
\end{split}
\end{equation*}
Furthermore, for $\|M_i\|_{\tilde{L}_{t}^{\infty}(\dot{B}^{\frac{N}{p}}_{p,1})}(i=1,2)$, by (i) we have
\begin{equation*}
\begin{split}&\|M_i\|_{\tilde{L}_{t}^{\infty}(\dot{B}^{\frac{N}{p}}_{p,1})}\\&\leq C\Big(1+\|(f_1,f_2)\|_{L_{t}^{\infty}(L^{\infty})}\|(g_1,g_2)\|_{L_{t}^{\infty}(L^{\infty})}\Big)^{[\frac{N}{p}]+1} \Big(\|(f_1,g_1)\|_{\tilde{L}_{t}^{\infty}(\dot{B}^{\frac{N}{p}}_{p,1})} +\|(f_2,g_2)\|_{\tilde{L}_{t}^{\infty}(\dot{B}^{\frac{N}{p}}_{p,1})}\Big).
\end{split}
\end{equation*}
Then, we conclude  \eqref{eq:2.3}.

We finish this subsection by listing an elementary but useful
inequality.
\begin{Lemma}\cite{MN2}\label{lemma2.13}\quad  Let $r_1,r_2>0$ satisfy $\max\{r_1,r_2\}>1$. Then
$$\int_0^t(1+t-\tau)^{-r_1}(1+\tau)^{-r_2}d\tau\leq C(r_1,r_2)(1+t)^{-\min\{r_1,r_2\}}.$$
\end{Lemma}
\section{Maximal Regularity Estimates for  the Linearized  System in Besov Spaces}
\ \ \ \ \  In this section,  we will exhibit the parabolic properties of the linearized system for all frequencies of the Cauchy problem  \eqref{equ:CTFS2}-\eqref{equ:CTFS3}.  The
key to these remarkable properties is given by the following lemma.
\begin{Lemma}\label{all part} Let $T\geq0$,
$s\in \mathbb{R}$, $1\leq r\leq q\leq\infty$.
Assume that $(c^{+},\, u^{+},\, c^{-},\, u^{-})$ is a solution to system \eqref{equ:CTFS2}-\eqref{equ:CTFS3} on $[0,T]\times \mathbb{R}^{N}$,  then,
\begin{equation}\label{regular estimate}\begin{split}
&\big\|\big((\sqrt{\beta_1}+\Lambda)c^{+},u^{+},(\sqrt{\beta_4}+\Lambda)c^{-},u^{-}\big)\big\|
_{\wt L^q_t(\dot B^{s+\frac{2}{q}}_{2,1})}
\\&\leq C \big\|\big((\sqrt{\beta_1}+\Lambda)c^{+}_0,u^{+}_0,(\sqrt{\beta_4}+\Lambda)c^{-}_0,u^{-}_0\big)\big\|
_{\dot{B}^{s}_{2,1}}\\&\quad+\big\|\big((\sqrt{\beta_1}+\Lambda)H_{1},H_{2},(\sqrt{\beta_4}+\Lambda)H_{3},H_{4}\big)
\big\|_{\tilde{L}_{t}^{r}(\dot{B}^{s-2-\frac{2}{r}}_{2,1})}  \quad \hbox{for} \quad \forall\, t\in[0,T].
\end{split}\end{equation}
\end{Lemma}

\noindent{\bf Proof.}\,   Applying the orthogonal projectors    $\cP$ and $\cQ$  over divergence-free and potential vector-fields, respectively,
 to the second equation and the fourth equation,
and setting  $\nu^{\pm}\eqdefa\nu_{1}^{\pm}+\nu_{2}^{\pm}$, System \eqref{equ:CTFS2} translates into
\begin{equation}\label{gama2}
\left\{
\begin{aligned}{}
&\p_t\cP u^{+}-\nu_{1}^{+}\Delta \cP u^{+}=\cP H_{2},
\\&\p_t\cP u^{-}-\nu_{1}^{-}\Delta \cP u^{-}=\cP H_{4},
\end{aligned}
\right.
\end{equation}
and
\begin{equation}\label{gama3}
\left\{
\begin{aligned}{}
&\p_tc^{+}+\textrm{div}\cQ u^{+}=H_{1},\\
&\p_t\cQ u^{+}+\beta_{1}\nabla c^{+}
+\beta_{2}\nabla c^{-}-\nu^{+}\Delta\cQ u^{+}-\nabla \Delta c^{+}=\cQ H_{2},
\\&\p_tc^{-}+\textrm{div}\cQ u^{-}=H_{3},
\\&\p_t\cQ u^{-}+\beta_{3}\nabla c^{+}
+\beta_{4}\nabla c^{-}-\nu^{-}\Delta \cQ u^{-}-\nabla \Delta c^{-}=\cQ H_{4}.
\end{aligned}
\right.
\end{equation}
Note that taking advantage of  Duhamel's formula reduces the proof to the case
where $H_{1}\equiv0$, $H_{2}\equiv0$, $H_{3}\equiv0$ and $H_{4}\equiv0$. Thus, from \eqref{gama2} we readily get, after taking the (space) Fourier transform,
\begin{equation}\label{gama4}
\frac{1}{2}\frac{d}{dt}\big(|\widehat{\cP u^{+}}|^{2}
+|\widehat{\cP u^{-}}|^{2}\big)
+\nu_{1}^{+}|\xi|^{2}|\widehat{\cP u^{+}}|^{2}
+\nu_{1}^{-}|\xi|^{2}|\widehat{\cP u^{-}}|^{2}=0.
\end{equation}
On the other hand, it
is convenient to set $d^{\pm}:=\Lambda^{-1}\textrm{div}u^{\pm}$ (with $\cF(\Lambda^s u)(\xi)\eqdefa |\xi|^s\wh u(\xi)$),
keeping in mind  that  bounding $d^{\pm}$ or $\cQ u^{\pm}$  is equivalent,
as one can go from $d^{\pm}$ to $\cQ u^{\pm}$ or from $\cQ u^{\pm}$ to $d^{\pm}$ by means of
a $0$ order homogeneous Fourier multiplier. Hence, from System \eqref{gama3}, we discover that $(c^{\pm}, d^{\pm})$ satisfies
\begin{equation}\label{gama5}
\left\{
\begin{aligned}{}
&\p_tc^{+}+\Lambda d^{+}=0,\\
&\p_td^{+}-\beta_{1}\Lambda c^{+}
-\beta_{2}\Lambda c^{-}-\nu^{+}\Delta d^{+}-\Lambda^{3} c^{+}=0,
\\&\p_tc^{-}+\Lambda d^{-}=0,
\\&\p_td^{-}-\beta_{3}\Lambda c^{+}
-\beta_{4}\Lambda c^{-}-\nu^{-}\Delta d^{-}-\Lambda^{3} c^{-}=0.
\end{aligned}
\right.
\end{equation}
By using the Fourier transform, from \eqref{gama5}, we have
\begin{equation}\label{gama6}
\left\{
\begin{aligned}{}
&\p_t\widehat{c^{+}}+|\xi|\widehat{d^{+}}=0,\\
&\p_t\widehat{d^{+}}-\beta_{1}|\xi| \widehat{c^{+}}
-\beta_{2}|\xi| \widehat{c^{-}}+\nu^{+}|\xi|^{2} \widehat{d^{+}}-|\xi|^{3} \widehat{c^{+}}=0,
\\&\p_t\widehat{c^{-}}+|\xi|\widehat{d^{-}}=0,
\\&\p_t\widehat{d^{-}}-\beta_{3}|\xi|\widehat{ c^{+}}
-\beta_{4}|\xi|\widehat{c^{-}}+\nu^{-}|\xi|^{2} \widehat{d^{-}}-|\xi|^{3} \widehat{c^{-}}=0.
\end{aligned}
\right.
\end{equation}
Applying the energy argument of Godunov \cite{GO}
 for partially dissipative first-order symmetric systems( and developed further in \cite{KF}) to the system \eqref{gama6}.
Multiplying the first equation in \eqref{gama6} by the conjugate $\overline{\wh c^{+}}$ of $\wh c^{+},$
 the second one by $\overline{\wh d^{+}}$, the third one by $\overline{\wh c^{-}}$ and the fourth one by $\overline{\wh d^{-}}$,
we get
\begin{equation}\label{R-E19}
\frac12\frac{d}{dt}|\wh c^{+}|^2+|\xi| \mathrm{Re}(\wh c^{+}\,\ov{\wh d^{+}})=0
\end{equation}
and, because  $\mathrm{Re}(\wh c^{+}\,\ov{\wh d^{+}})=\mathrm{Re}(\ov{\wh c^{+}}\,\wh d^{+}),$
\begin{equation}\label{R-E20}
\frac12\frac{d}{dt}|\wh d^{+}|^2+\nu^{+}|\xi|^2|\wh d^{+}|^2-|\xi|(\beta_1+|\xi|^2)\mathrm{Re}(\wh c^{+}\,\ov{\wh d^{+}})-\beta_2|\xi|\mathrm{Re}(\wh c^{-}\,\ov{\wh d^{+}})=0,
\end{equation}
\begin{equation}\label{R-E21}
\frac12\frac{d}{dt}|\wh c^{-}|^2+|\xi| \mathrm{Re}(\wh c^{-}\,\ov{\wh d^{-}})=0,
\end{equation}
and
\begin{equation}\label{R-E22}
\frac12\frac{d}{dt}|\wh d^{-}|^2+\nu^{-}|\xi|^2|\wh d^{-}|^2-|\xi|(\beta_4+|\xi|^2)\mathrm{Re}(\wh c^{-}\,\ov{\wh d^{-}})-\beta_3|\xi|\mathrm{Re}(\wh c^{+}\,\ov{\wh d^{-}})=0.
\end{equation}
Further, multiplying the first and third
equations of \eqref{gama6} by $\overline{\wh c^{-}}$ and $\overline{\wh c^{+}}$ respectively, we have
\begin{equation}\label{R-E23}
\frac{d}{dt}\mathrm{Re}(\wh c^{+}\,\ov{\wh c^{-}})+|\xi| \mathrm{Re}(\wh c^{-}\,\ov{\wh d^{+}})+|\xi| \mathrm{Re}(\wh c^{+}\,\ov{\wh d^{-}})=0.
\end{equation}
Adding \eqref{R-E19}$\times(\beta_1+|\xi|^2)$, \eqref{R-E21}$\times(\beta_4+|\xi|^2)$ and  \eqref{R-E23}$\times \beta_2$ to \eqref{R-E20} and \eqref{R-E22},  and using $\beta_2=\beta_3$, we obtain
\begin{equation}\label{R-E24}
\begin{aligned}
\frac12\frac{d}{dt}\Big((\beta_1&+|\xi|^2)|\wh c^{+}|^2+(\beta_4+|\xi|^2)|\wh c^{-}|^2+2\beta_2\mathrm{Re}(\wh c^{+}\,\ov{\wh c^{-}}) +|\wh d^{+}|^2+|\wh d^{-}|^2\Big)
\\&+\nu^{+}|\xi|^2|\wh d^{+}|^2+\nu^{-}|\xi|^2|\wh d^{-}|^2 =0.
\end{aligned}
\end{equation}
In order to track  dissipations arising for $c^{+}$ and $c^{-}$ , let us multiply the first and second
equations of \eqref{gama6} by $-|\xi|\ov{\wh d^{+}}$ and  $-|\xi|\ov{\wh c^{+}}$, the third one by  $-|\xi|\ov{\wh d^{-}}$ and the fourth one by $\overline{\wh c^{-}}$, respectively. Adding them, we get
 \begin{equation}\label{R-E25}
\begin{aligned}
&\frac{d}{dt}\Big(-|\xi|\mathrm{Re}(\wh c^{+}\,\ov{\wh d^{+}})-|\xi|\mathrm{Re}(\wh c^{+}\,\ov{\wh d^{+}})\Big)
+|\xi|^2(\beta_1+|\xi|^2)|\wh c^{+}|^2+|\xi|^2(\beta_4+|\xi|^2)|\wh c^{-}|^2\\&\quad
-\beta_2|\xi|^2\mathrm{Re}(\wh c^{+}\,\ov{\wh c^{-}})-\beta_3|\xi|^2\mathrm{Re}(\wh c^{+}\,\ov{\wh c^{-}})-\nu^{+}|\xi|^3\mathrm{Re}(\wh c^{+}\,\ov{\wh d^{+}})-\nu^{-}|\xi|^3\mathrm{Re}(\wh c^{-}\,\ov{\wh d^{-}})
\\&\quad-|\xi|^2|\wh d^{+}|^2-|\xi|^2|\wh d^{-}|^2 =0.
\end{aligned}
\end{equation}
Adding \eqref{R-E19}$\times\nu^{+}|\xi|^2$  and \eqref{R-E21}$\times\nu^{-}|\xi|^2$ to \eqref{R-E25}, we get
 \begin{equation}\label{R-E26}
\begin{aligned}
&\frac{1}{2}\frac{d}{dt}\Big(\nu^{+}|\xi|^2|\wh c^{+}|^2 +\nu^{-}|\xi|^2 |\wh c^{-}|^2-2|\xi|\mathrm{Re}(\wh c^{+}\,\ov{\wh d^{+}})-2|\xi|\mathrm{Re}(\wh c^{+}\,\ov{\wh d^{+}})\Big)
\\&\quad+|\xi|^2(\beta_1+|\xi|^2)|\wh c^{+}|^2+|\xi|^2(\beta_4+|\xi|^2)|\wh c^{-}|^2
-\beta_2\mathrm{Re}(\wh c^{+}\,\ov{\wh c^{-}})-\beta_3\mathrm{Re}(\wh c^{+}\,\ov{\wh c^{-}})
\\&\quad-|\xi|^2|\wh d^{+}|^2-|\xi|^2|\wh d^{-}|^2 =0.
\end{aligned}
\end{equation}
Therefore, by multiplying \eqref{R-E26} by a small enough constant $\delta> 0$ (to be determined
later) and adding it to \eqref{R-E24}, we get
 \begin{equation}\label{R-E27}
\begin{aligned}
&\frac{1}{2}\frac{d}{dt}\Big((\beta_1+|\xi|^2)|\wh c^{+}|^2+(\beta_4+|\xi|^2)|\wh c^{-}|^2+2\beta_2\mathrm{Re}(\wh c^{+}\,\ov{\wh c^{-}}) +|\wh d^{+}|^2+|\wh d^{-}|^2\\&\quad+\delta\nu^{+}|\xi|^2|\wh c^{+}|^2 +\delta\nu^{-}|\xi|^2 |\wh c^{-}|^2-2\delta|\xi|\mathrm{Re}(\wh c^{+}\,\ov{\wh d^{+}})-2\delta|\xi|\mathrm{Re}(\wh c^{-}\,\ov{\wh d^{-}})\Big)
\\&\quad+\delta|\xi|^2(\beta_1+\delta|\xi|^2)|\wh c^{+}|^2+\delta|\xi|^2(\beta_4+|\xi|^2)|\wh c^{-}|^2\\&\quad
-\delta\beta_2|\xi|^2\mathrm{Re}(\wh c^{+}\,\ov{\wh c^{-}})-\delta\beta_3|\xi|^2\mathrm{Re}(\wh c^{+}\,\ov{\wh c^{-}})
\\&\quad+(\nu^{+}-\delta)|\xi|^2|\wh d^{+}|^2+(\nu^{+}-\delta)|\xi|^2|\wh d^{-}|^2=0.
\end{aligned}
\end{equation}
By Young's inequality, for $\delta<\frac{1}{3}$,   we have
 \begin{equation*}\begin{aligned}\Big|2\delta|\xi|\mathrm{Re}(\wh c^{+}\,\ov{\wh d^{+}})\Big|&\leq2\delta\big(\frac{|\xi|^2|\wh c^{+}|^2}{2}+ \frac{1}{2} |\wh d^{+}|^2\big)
\\&\leq\frac{1}{3}|\xi|^2|\wh c^{+}|^2+ \frac{1}{3} |\wh d^{+}|^2,\end{aligned}\end{equation*}
 \begin{equation*}\begin{aligned}\Big|2\delta|\xi|\mathrm{Re}(\wh c^{-}\,\ov{\wh d^{-}})\Big|&\leq2\delta\big(\frac{|\xi|^2|\wh c^{-}|^2}{2}+ \frac{1}{2} |\wh d^{-}|^2\big)
\\&\leq\frac{1}{3}|\xi|^2|\wh c^{-}|^2+ \frac{1}{3} |\wh d^{-}|^2.\end{aligned}\end{equation*}
Using further $\beta^{2}_2=\beta^{2}_3=\beta_{1}\beta_{4}$ and Young's inequality,  and choosing $M=\frac{\beta_1}{\beta_2}$,  we get
\begin{equation*}\begin{aligned}\Big|2\beta_2\mathrm{Re}(\wh c^{+}\,\ov{\wh c^{-}})\Big|&\leq M\beta_2|\wh c^{+}|^2+ \frac{\beta_2}{M} |\wh c^{-}|^2
\\&\leq \beta_1|\wh c^{+}|^2+ \beta_4|\wh c^{-}|^2.\end{aligned}\end{equation*}
Thus,
$$\mathcal{L}^{2}\approx (\beta_1+|\xi|^2)|\wh c^{+}|^2+(\beta_4+|\xi|^2)|\wh c^{-}|^2+|\wh d^{+}|^2+|\wh d^{-}|^2,$$
where $\mathcal{L}^{2}:=\mathcal{L}^{2}(\xi,t)$ is a Lyapunov functional defined by
\begin{equation*}
\begin{split}
\mathcal{L}^{2}=&(\beta_1+|\xi|^2)|\wh c^{+}|^2+(\beta_4+|\xi|^2)|\wh c^{-}|^2+2\beta_2\mathrm{Re}(\wh c^{+}\,\ov{\wh c^{-}}) +|\wh d^{+}|^2+|\wh d^{-}|^2\\&\quad+\delta\nu^{+}|\xi|^2|\wh c^{+}|^2 +\delta\nu^{-}|\xi|^2 |\wh c^{-}|^2-2\delta|\xi|\mathrm{Re}(\wh c^{+}\,\ov{\wh d^{+}})-2\delta|\xi|\mathrm{Re}(\wh c^{-}\,\ov{\wh d^{-}}).
\end{split}
\end{equation*}
On the other hand, taking $\delta=\min\{\frac{1}{3},\nu^{+},\nu^{-} \}$, there exists a positive constant $C$ such that
\begin{equation*}\label{R-E27}
\begin{aligned}
&\delta|\xi|^2(\beta_1+\delta|\xi|^2)|\wh c^{+}|^2+\delta|\xi|^2(\beta_4+|\xi|^2)|\wh c^{-}|^2\\&\qquad
-\delta\beta_2|\xi|^2\mathrm{Re}(\wh c^{+}\,\ov{\wh c^{-}})-\delta\beta_3|\xi|^2\mathrm{Re}(\wh c^{+}\,\ov{\wh c^{-}})
\\&\qquad+(\nu^{+}-\delta)|\xi|^2|\wh d^{+}|^2+(\nu^{+}-\delta)|\xi|^2|\wh d^{-}|^2
\\&\quad\geq C|\xi|^2\Big((\beta_1+|\xi|^2)|\wh c^{+}|^2+(\beta_4+|\xi|^2)|\wh c^{-}|^2+|\wh d^{+}|^2+|\wh d^{-}|^2\Big).
\end{aligned}
\end{equation*}
Thus,
\begin{equation}\label{gama7}
\frac{d}{dt}\mathcal{L}^{2}+C|\xi|^{2}\mathcal{L}^{2}\leq0.
\end{equation}
Combining \eqref{gama4} and  \eqref{gama7}, we get
\begin{equation}\label{gama8}
\big|(\sqrt{\beta_1}\widehat{c^{+}},|\xi|\widehat{c^{+}}, \widehat{u^{+}}, \sqrt{\beta_4}\widehat{c^{-}},|\xi|\widehat{c^{-}}, \widehat{u^{-}})\big|
\leq Ce^{-c_{0}|\xi|^{2}t}\big|(\sqrt{\beta_1}\widehat{c^{+}},|\xi|\widehat{c^{+}}, \widehat{u^{+}}, \sqrt{\beta_4}\widehat{c^{-}},|\xi|\widehat{c^{-}}, \widehat{u^{-}})\big|(0).
\end{equation}
Granted with the above estimates, by means of the definition of  Littlewood-Paley decomposition and  Fourier-Plancherel
theorem, we get
\begin{equation}\label{gama9}\begin{split}
&\big\|\big((\sqrt{\beta_1}+\Lambda)c^{+},u^{+},(\sqrt{\beta_4}+\Lambda)c^{-},u^{-}\big)\big\|
_{\wt L^r_t(\dot B^{s+\frac{2}{r}}_{2,1})}
\\&\quad\leq C\big\|\big((\sqrt{\beta_1}+\Lambda)c^{+}_0,u^{+}_0,(\sqrt{\beta_4}+\Lambda)c^{-}_0,u^{-}_0\big)\big\|
_{\dot{B}^{s}_{2,1}}.
\end{split}\end{equation}
Then, for general source terms,  from Duhamel's formula we finally obtain \eqref{regular estimate}.
This completes the proof of Lemma \ref{all part}.
\section{The Unique Global Solvability}

\subsection{Global \emph{a priori} estimates}
The subsection is devoted to  exploiting an important global \emph{a priori}  estimates for the Cauchy problem \eqref{equ:CTFS2}-\eqref{equ:CTFS3}.
\begin{Proposition}\label{ab} Let $T\geq0$,
$N\geq2$, and  $(c^{+},\,u^{+},\,c^{-},\,u^{-})$ be a solution to the Cauchy  problem \eqref{equ:CTFS2}-\eqref{equ:CTFS3} on $[0,T]\times \mathbb{R}^{N}$,
we have
\begin{equation}\label{eq:5.1}X(t)\leq C\Big(X(0)+ \big(1+X^{2}(t)\big)^{[\frac{N}{2}]+1}\big(X^{2}(t)+X^{3}(t)\big)\Big), \quad \hbox{for} \quad \forall \,t\in[0,T],\end{equation}
where \begin{equation*}\label{eq:5.2}\begin{split}
X(t)&\eqdefa\big\|\big((\sqrt{\beta_1}+\Lambda)c^{+},u^{+},(\sqrt{\beta_4}+\Lambda)c^{-},u^{-}\big)\big\|
_{\wt L^\infty_t(\dot B^{\frac{N}{2}-1}_{2,1})}\\&\qquad+\big\|\big((\sqrt{\beta_1}+\Lambda)c^{+},u^{+},(\sqrt{\beta_4}+\Lambda)c^{-},u^{-}\big)\big\|
_{\wt L^1_t(\dot B^{\frac{N}{2}+1}_{2,1})}.
\end{split}\end{equation*}
\end{Proposition}
\noindent{\bf Proof.}\, In Lemma \ref{all part},  taking $q=1,\infty$, $r=1$ and $s=\frac{N}{2}-1$,  for  $t\in[0,T]$,   we have
\begin{equation}\label{eq:5.2}\begin{split}
X(t)&\leq C \big\|\big((\sqrt{\beta_1}+\Lambda)c^{+}_0,u^{+}_0,(\sqrt{\beta_4}+\Lambda)c^{-}_0,u^{-}_0\big)\big\|
_{\dot{B}^{\frac{N}{2}-1}_{2,1}}\\&\qquad+\big\|\big((\sqrt{\beta_1}+\Lambda)H_{1},H_{2},(\sqrt{\beta_4}+\Lambda)H_{3},H_{4}\big)
\big\|_{L_{t}^{1}(\dot{B}^{\frac{N}{2}-1}_{2,1})}.
\end{split}\end{equation}
In what follows, we derive some nonlinear estimates from the terms $\|((\sqrt{\beta_1}+\Lambda)H_{1},H_{2},(\sqrt{\beta_4}+\Lambda)H_{3},H_{4})
\|_{L_{t}^{1}(\dot{B}^{\frac{N}{2}-1}_{2,1})}$. First, by Proposition \ref{p26} and the embedding $\dot{B}^\frac N2_{2,1}\hookrightarrow
L^\infty$,  we have
\begin{equation}\label{eq:5.3}
\begin{split}&\|(\sqrt{\beta_1}+\Lambda)H_{1}\|_{L_{t}^{1}(\dot{B}^{\frac{N}{2}-1}_{2,1})}\\
&\quad\leq C\|(\sqrt{\beta_1}+\Lambda)(c^{+}u^{+})\|_{L_{t}^{1}(\dot{B}^{\frac{N}{2}}_{2,1})}
\\&\quad\leq C\|(c^{+}u^{+})\|_{L_{t}^{1}(\dot{B}^{\frac{N}{2}}_{2,1})}+C\|(c^{+}u^{+})\|_{L_{t}^{1}(\dot{B}^{\frac{N}{2}+1}_{2,1})}
\\&\quad\leq C\|c^{+}\|_{L_{t}^{2}(\dot{B}^{\frac{N}{2}}_{2,1})}\|u^{+}\|_{L_{t}^{2}(\dot{B}^{\frac{N}{2}}_{2,1})}
+C\|c^{+}\|_{L_{t}^{\infty}(L^{\infty})}\|u^{+}\|_{L_{t}^{1}(\dot{B}^{\frac{N}{2}+1}_{2,1})}
\\&\qquad+C\|c^{+}\|_{L_{t}^{2}(\dot{B}^{\frac{N}{2}+1}_{2,1})}\|u^{+}\|_{L_{t}^{2}(L^{\infty})}
\\&\quad\leq C\|c^{+}\|_{L_{t}^{2}(\dot{B}^{\frac{N}{2}}_{2,1})}\|u^{+}\|_{L_{t}^{2}(\dot{B}^{\frac{N}{2}}_{2,1})}
+C\|c^{+}\|_{L_{t}^{\infty}(\dot{B}^{\frac{N}{2}}_{2,1})}\|u^{+}\|_{L_{t}^{1}(\dot{B}^{\frac{N}{2}+1}_{2,1})}
\\&\quad\quad+C\|c^{+}\|_{L_{t}^{2}(\dot{B}^{\frac{N}{2}+1}_{2,1})}\|u^{+}\|_{L_{t}^{2}(\dot{B}^{\frac{N}{2}}_{2,1})}
\\&\quad\leq CX^{2}(t),
\end{split}
\end{equation}
where we have used the following interpolation inequalities,
$$\|c\|_{\tilde{L}^2_t(\dot B^{\frac N2}_{2,1})}\leq \Big(\|c\|_{\tilde{L}^\infty_t(\dot B^{\frac N2-1}_{2,1})}\Big)^{\frac{1}{2} }\Big(\|c\|_{L^1_t(\dot B^{\frac N2+1}_{2,1})}\Big)^{\frac{1}{2}}, $$
$$\|c\|_{\tilde{L}^2_t(\dot B^{\frac N2+1}_{2,1})}\leq \Big(\|c\|_{\tilde{L}^\infty_t(\dot B^{\frac N2}_{2,1})}\Big)^{\frac{1}{2} }\Big(\|c\|_{L^1_t(\dot B^{\frac N2+2}_{2,1})}\Big)^{\frac{1}{2}}, $$
and
$$\|u\|_{\tilde{L}^2_t(\dot B^{\frac N2}_{2,1})}\leq \Big(\|u\|_{\tilde{L}^\infty_t(\dot B^{\frac N2-1}_{2,1})}\Big)^{\frac{1}{2} }\Big(\|u\|_{L^1_t(\dot B^{\frac N2+1}_{2,1})}\Big)^{\frac{1}{2}}.$$
Similarly,
\begin{equation}\label{eq:5.4}
\begin{split}\|(\sqrt{\beta_4}+\Lambda)H_{3}\|_{L_{t}^{1}(\dot{B}^{\frac{N}{2}-1}_{2,1})}\leq CX^{2}(t).
\end{split}
\end{equation}
Next, we bound the term $H_{2}$.
By Lemma \ref{p27} (i) , we get
\begin{equation*}\label{318}
\begin{split}\big\|g_{+}(c^{+},c^{-})\big\|_{\tilde{L}_{t}^{2}(\dot{B}^\frac N2_{2,1})}
\leq C\Big(1+\|c^{+}\|_{L_{t}^\infty(L^\infty)}\|c^{-}\|_{L_{t}^\infty(L^\infty)}\Big)^{[\frac{N}{2}]+1}\|(c^{+},c^{-)}\|_{\tilde{L}_{t}^{2}(\dot{B}^\frac N2_{2,1})},
\end{split}
\end{equation*}
\begin{equation*}\label{318}
\begin{split}\big\|\tilde{g}_{+}(c^{+},c^{-})\big\|_{\tilde{L}_{t}^{2}(\dot{B}^\frac N2_{2,1})}
\leq C\Big(1+\|c^{+}\|_{L_{t}^\infty(L^\infty)}\|c^{-}\|_{L_{t}^\infty(L^\infty)}\Big)^{[\frac{N}{2}]+1}\|(c^{+},c^{-)}\|_{\tilde{L}_{t}^{2}(\dot{B}^\frac N2_{2,1})},
\end{split}
\end{equation*}

\begin{equation*}\label{318}
\begin{split}\big\|l_{+}(c^{+},c^{-})\big\|_{\tilde{L}_{t}^{\infty}(\dot{B}^\frac N2_{2,1})}
\leq C\Big(1+\|c^{+}\|_{L_{t}^\infty(L^\infty)}\|c^{-}\|_{L_{t}^\infty(L^\infty)}\Big)^{[\frac{N}{2}]+1}\|(c^{+},c^{-)}\|_{\tilde{L}_{t}^{\infty}(\dot{B}^\frac N2_{2,1})},
\end{split}
\end{equation*}
\begin{equation*}\label{318}
\begin{split}&\Big\|h_{+}(c^{+},c^{-})-\frac{(\mathcal{C}^{2}\alpha^{-})(1,1)}
{s_{-}^{2}(1,1)}\Big\|_{\tilde{L}_{t}^{\infty}(\dot{B}^\frac N2_{2,1})}
\\&\quad\leq C\Big(1+\|c^{+}\|_{L_{t}^\infty(L^\infty)}\|c^{-}\|_{L_{t}^\infty(L^\infty)}\Big)^{[\frac{N}{2}]+1}\|(c^{+},c^{-)}\|_{\tilde{L}_{t}^{\infty}(\dot{B}^\frac N2_{2,1})},
\end{split}
\end{equation*}
and
\begin{equation*}\label{318}
\begin{split}&\Big\|k_{+}(c^{+},c^{-})-\frac{\mathcal{C}^{2}(1,1)}
{s_{+}^{2}\rho^{+}(1,1)}\Big\|_{\tilde{L}_{t}^{\infty}(\dot{B}^\frac N2_{2,1})}
\\&\quad\leq C\Big(1+\|c^{+}\|_{L_{t}^\infty(L^\infty)}\|c^{-}\|_{L_{t}^\infty(L^\infty)}\Big)^{[\frac{N}{2}]+1}\|(c^{+},c^{-)}\|_{\tilde{L}_{t}^{\infty}(\dot{B}^\frac N2_{2,1})}.
\end{split}
\end{equation*}
Based on these bounds, Thanks to Proposition \ref{p26}, Lemma \ref{p27} (i) and   the embedding $\dot{B}^\frac N2_{2,1}\hookrightarrow
L^\infty$, we easily infer the following estimates
\begin{equation*}\label{318}
\begin{split}&\big\|g_{+}(c^{+},c^{-})\partial_{i}c^{+}
-\tilde{g}_{+}(c^{+},c^{-})\partial_{i}c^{-}\big\|_{L^1_{t}(\dot{B}^{\f{N}{2}-1}_{2,1})}
\\&\quad\leq C\big\|g_{+}(c^{+},c^{-})\big\|_{L^2_{t}(\dot{B}^{\frac{N}{2}}_{2,1})}
\big\|\partial_{i}c^{+}\big\|_{L^2_{t}(\dot{B}^{\f{N}{2}-1}_{2,1})}+\big\|\tilde{g}_{+}(c^{+},c^{-})\big\|_{L^2_{t}(\dot{B}^{\frac{N}{2}}_{2,1})}
\big\|\partial_{i}c^{-}\big\|_{L^2_{t}(\dot{B}^{\f{N}{2}-1}_{2,1})}
\\&\quad\leq C\big\|g_{+}(c^{+},c^{-})\big\|_{L^2_{t}(\dot{B}^{\frac{N}{2}}_{2,1})}
\|c^{+}\|_{L^2_{t}(\dot{B}^{\f{N}{2}}_{2,1})}+\|\tilde{g}_{+}(c^{+},c^{-})\|_{L^2_{t}(\dot{B}^{\frac{N}{2}}_{2,1})}
\|c^{-}\|_{L^2_{t}(\dot{B}^{\f{N}{2}}_{2,1})}
\\&\quad\leq C\Big(1+\|c^{+}\|_{L_{t}^\infty(L^\infty)}\|c^{-}\|_{L_{t}^\infty(L^\infty)}\Big)^{[\frac{N}{2}]+1}\|(c^{+},c^{-)}\|_{\tilde{L}_{t}^{2}(\dot{B}^\frac N2_{2,1})}^{2}
\\&\quad\leq C \Big(1+X^{2}(t)\Big)^{[\frac{N}{2}]+1}X^{2}(t),
\end{split}
\end{equation*}
\begin{equation*}\label{318}
\begin{split}\big\|u^{+}\cdot\nabla u_{i}^{+}\big\|_{L^1_{T_0}(\dot{B}^{\f{N}{2}-1}_{2,1})}
&\leq C \|u^{+}\|_{L^\infty_{t}(\dot{B}^{\f{N}{2}-1}_{2,1})}\|\nabla u^{+}\|_{L^1_{t}(\dot{B}^{\f{N}{2}}_{2,1})}
\\&\leq C \|u^{+}\|_{L^\infty_{t}(\dot{B}^{\f{N}{2}-1}_{2,1})}\| u^{+}\|_{L^1_{t}(\dot{B}^{\f{N}{2}+1}_{2,1})}
\\&\leq C X^{2}(t),
\end{split}
\end{equation*}
\begin{equation*}\label{318}
\begin{split}&\big\|l_{+}(c^{+},c^{-})\partial_{j}^{2}u_{i}^{+}\big\|_{L^1_{t}(\dot{B}^{\f{N}{2}-1}_{2,1})}
\\&\quad\leq C\big \|l_{+}(c^{+},c^{-})\big\|_{\tilde{L}_{t}^{\infty}(\dot{B}^\frac N2_{2,1})}
\big\|\partial_{j}^{2}u_{i}^{+}\big\|_{L^1_{t}(\dot{B}^{\f{N}{2}-1}_{2,1})}
\\&\quad\leq C\Big(1+\|c^{+}\|_{L_{t}^\infty(L^\infty)}\|c^{-}\|_{L_{t}^\infty(L^\infty)}\Big)^{[\frac{N}{2}]+1}\|(c^{+},c^{-)}\|_{\tilde{L}_{t}^{\infty}(\dot{B}^\frac N2_{2,1})}
\|u^{+}\|_{L^1_{t}(\dot{B}^{\f{N}{2}+1}_{2,1})}
\\&\quad\leq C \Big(1+X^{2}(t)\Big)^{[\frac{N}{2}]+1}X^{2}(t),
\end{split}
\end{equation*}
\begin{equation*}\label{318}
\begin{split}&\big\|h_{+}(c^{+},c^{-})\partial_{j}c^{+}\partial_{j}u^{+}_{i}\big\|_{L^1_{t}(\dot{B}^{\f{N}{2}-1}_{2,1})}
\\&\quad\leq C \big\|h_{+}(c^{+},c^{-})\big\|_{\tilde{L}_{t}^{\infty}(\dot{B}^\frac N2_{2,1})}
\big\|\partial_{j}c^{+}\partial_{j}u^{+}_{i}\big\|_{L^1_{t}(\dot{B}^{\f{N}{2}-1}_{2,1})}
\\&\quad\leq C\Big(1+\|c^{+}\|_{L_{t}^\infty(L^\infty)}\|c^{-}\|_{L_{t}^\infty(L^\infty)}\Big)^{[\frac{N}{2}]+1}\|(c^{+},c^{-)}\|_{\tilde{L}_{t}^{\infty}(\dot{B}^\frac N2_{2,1})}
\|\nabla c^{+}\|_{L^\infty_{t}(\dot{B}^{\f{N}{2}-1}_{2,1})}\|\nabla u^{+}\|_{L^1_{t}(\dot{B}^{\f{N}{2}}_{2,1})}
\\&\quad\leq C \Big(1+X^{2}(t)\Big)^{[\frac{N}{2}]+1}X^{3}(t),
\end{split}
\end{equation*} and
\begin{equation*}\label{318}
\begin{split}&\Big\|k_{+}(c^{+},c^{-})\partial_{j}c^{-}\partial_{j}u^{+}_{i}\Big\|_{L^1_{t}(\dot{B}^{\f{N}{2}-1}_{2,1})}
\\&\quad\leq C\Big(1+\|k_{+}(c^{+},c^{-})-\frac{\mathcal{C}^{2}(1,1)}
{s_{+}^{2}\rho^{+}(1,1)}\|_{\tilde{L}_{t}^{\infty}(\dot{B}^\frac N2_{2,1})}\Big)
\big\|\partial_{j}c^{+}\partial_{j}u^{+}_{i}\big\|_{L^1_{t}(\dot{B}^{\f{N}{2}-1}_{2,1})}
\\&\quad\leq C\Big(1+\|c^{+}\|_{L_{t}^\infty(L^\infty)}\|c^{-}\|_{L_{t}^\infty(L^\infty)}\Big)^{[\frac{N}{2}]+1}\|(c^{+},c^{-)}\|_{\tilde{L}_{t}^{\infty}(\dot{B}^\frac N2_{2,1})}
\|\nabla c^{-}\|_{L^\infty_{t}(\dot{B}^{\f{N}{2}-1}_{2,1})}\|\nabla u^{+}\|_{L^1_{t}(\dot{B}^{\f{N}{2}}_{2,1})}
\\&\quad\leq C \Big(1+X^{2}(t)\Big)^{[\frac{N}{2}]+1}X^{3}(t).
\end{split}
\end{equation*}
Hence, we gather that
\begin{equation}\label{eq:5.5}
\|H_{2}\|_{L^1_{t}(\dot{B}^{\f{N}{2}-1}_{2,1})}\le
C\big(1+X^{2}(t)\big)^{[\frac{N}{2}]+1}\big(X^{2}(t)+X^{3}(t)\big).
\end{equation}
Similarly,  we also have
\begin{equation}\label{eq:5.6}
\|H_{4}\|_{L^1_{t}(\dot{B}^{\f{N}{2}-1}_{2,1})}\le
C\big(1+X^{2}(t)\big)^{[\frac{N}{2}]+1}\big(X^{2}(t)+X^{3}(t)\big).
\end{equation}
 Substituting
\eqref{eq:5.3}-\eqref{eq:5.6} into  \eqref{eq:5.2}, we finally get  \eqref{eq:5.1}. This completes the proof of Proposition \ref{ab}.
\subsection{Global existence and uniqueness}
In order to solve the system \eqref{equ:CTFS2}-\eqref{equ:CTFS3} by fixed point theorem,  we define  the following  map
\begin{equation}\label{map}
\Phi: (c^{+},u^{+},c^{-},u^{-})\rightarrow (b^{+},v^{+},b^{-},v^{-})
\end{equation}
with $(b^{+},v^{+},b^{-},v^{-})$ the solution to
\begin{equation}\label{equ:linearCTFS2}
\left\{
\begin{aligned}{}
&\p_tb^{+}+\textrm{div}v^{+}=H_{1}(c^{+},u^{+}),\\
&\p_t{v}^{+}+\beta_{1}\nabla b^{+}
+\beta_{2}\nabla b^{-}-\nu_{1}^{+}\Delta v^{+}
-\nu_{2}^{+}\nabla\textrm{div}v^{+}-\nabla \Delta b^{+}=H_{2}(c^{+},u^{+},c^{-}),
\\&\p_tb^{-}+\textrm{div}v^{-}=H_{3}(c^{-},u^{-}),
\\&\p_tv^{-}+\beta_{3}\nabla b^{+}
+\beta_{4}\nabla b^{-}-\nu_{1}^{-}\Delta v^{-}
-\nu_{2}^{-}\nabla\textrm{div}v^{-}-\nabla \Delta b^{-}=H_{4}(c^{+},u^{-},c^{-}).
\end{aligned}
\right.
\end{equation}
Obviously, to prove the existence  part of the theorem,  we just have to show that $\Phi$ is a contraction map in a ball of $\mathbb{X}(t)$, where
$\mathbb{X}(t)\eqdefa\wt L^\infty_t(\dot B^{\frac{N}{2}-1}_{2,1})\cap  L^1_t(\dot B^{\frac{N}{2}+1}_{2,1})$ equipped with a norm
\begin{equation*}\label{991}\begin{split}\big\|(c^{+},u^{+},c^{-},u^{-})\big\|_{\mathbb{X}(t)}={X}(t).\end{split}\end{equation*}
We define a ball $B(O,R)$ centered at the origin by
\begin{equation}\label{ball}
B(0,R)=\Big\{(c^{+},u^{+},c^{-},u^{-})\in \mathbb{X}(t): \big\|(c^{+},u^{+},c^{-},u^{-})\big\|_{\mathbb{X}(t)}\leq R \Big\}.
\end{equation}
Assuming $R\leq1$, from Proposition \ref{ab} we have
\begin{equation}\label{991}\begin{split}&\big\|\Phi(c^{+},u^{+},c^{-},u^{-})\big\|_{\mathbb{X}(t)}\\&\quad\leq C \Big(X(0)+ \big(1+\big\|(c^{+},u^{+},c^{-},u^{-})\big\|^{2}_{\mathbb{X}(t)}\big)^{[\frac{N}{2}]+1}\\&\qquad\times\big(\big\|(c^{+},u^{+},c^{-},u^{-})\big\|^{2}_{\mathbb{X}(t)}
+\big\|(c^{+},u^{+},c^{-},u^{-})\big\|^{3}_{\mathbb{X}(t)}\big)\Big)
\\&\quad\leq C \Big(\eta+\big(1+R^{2}\big)^{[\frac{N}{2}]+1}\big(R^{2}+R^{3}\big)\Big)
\\&\quad\leq C \big(\eta+2R^{2}\big).
\end{split}\end{equation}
Choosing $(R,\eta)$ such that \begin{equation}\label{RRR}R\leq \min\{1, (4C)^{-1}\}\quad \hbox{and} \quad \eta\leq 2R^{2}.\end{equation}
Thus, from \eqref{991}, we finally deduce that $$\Phi(B(0,R))\subseteq B(0,R).$$
In order to show $\Phi$ is a  contraction map, one chooses two elements $(c^{+}_1,u^{+}_1,c^{-}_1,u^{-}_1)$ and $(c^{+}_2,u^{+}_2,c^{-}_2,u^{-}_2)$ in $B(0,R)$. According to \eqref{eq:5.2}, \eqref{map} and \eqref{equ:linearCTFS2},  we have
\begin{equation}\label{9911}\begin{split}&\big\|\Phi(c^{+}_1,u^{+}_1,c^{-}_1,u^{-}_1)-\Phi(c^{+}_2,u^{+}_2,c^{-}_2,u^{-}_2)\big\|_{\mathbb{X}(t)}
\\&\quad\leq C\Big(\big\|(\sqrt{\beta_1}+\Lambda)(H_{1}(c^{+}_1,u^{+}_1)-H_{1}(c^{+}_1,u^{+}_1))\big\|_{L_{t}^{1}(\dot{B}^{\frac{N}{2}-1}_{2,1})}
\\&\quad+\big\|H_{2}(c^{+}_1,u^{+}_1,c^{-}_1)-H_{2}(c^{+}_2,u^{+}_2,c^{-}_2)
\big\|_{L_{t}^{1}(\dot{B}^{\frac{N}{2}-1}_{2,1})}
\\&\quad+
\big\|\sqrt{\beta_4}+\Lambda)(H_{3}(c^{-}_1,u^{-}_1)-H_{3}(c^{-}_2,u^{-}_2))\big\|_{L_{t}^{1}(\dot{B}^{\frac{N}{2}-1}_{2,1})}
\\&\quad+
\big\|H_{4}(c^{+}_1,u^{-}_1,c^{-}_1)-H_{4}(c^{+}_2,u^{-}_2,c^{-}_2)
\big\|_{L_{t}^{1}(\dot{B}^{\frac{N}{2}-1}_{2,1})}\Big).
\end{split}\end{equation}
Further, using Lemma \ref{p27} (ii), similar to the estimate \eqref{eq:5.1}, we have
\begin{equation}\label{99111}\begin{split}&\big\|\Phi(c^{+}_1,u^{+}_1,c^{-}_1,u^{-}_1)-\Phi(c^{+}_2,u^{+}_2,c^{-}_2,u^{-}_2)\big\|_{\mathbb{X}(t)}
\\&\quad\leq C \Big(\big\|(c^{+}_1,u^{+}_1,c^{-}_1,u^{-}_1)\big\|_{\mathbb{X}(t)}
+\big\|(c^{+}_2,u^{+}_2,c^{-}_2,u^{-}_2)\big\|_{\mathbb{X}(t)}\Big)
\\&\qquad\times
\big\|(c^{+}_1-c^{+}_2,u^{+}_1-u^{+}_2,c^{-}_1-c^{-}_2,u^{-}_1-u^{-}_2)\big\|_{\mathbb{X}(t)}.
\end{split}\end{equation}
From  \eqref{RRR} we finally deduce that
\begin{equation}\label{99111}\big\|\Phi(c^{+}_1,u^{+}_1,c^{-}_1,u^{-}_1)-\Phi(c^{+}_2,u^{+}_2,c^{-}_2,u^{-}_2)\big\|_{\mathbb{X}(t)}\leq \frac{1}{2} \big\|(c^{+}_1-c^{+}_2,u^{+}_1-u^{+}_2,c^{-}_1-c^{-}_2,u^{-}_1-u^{-}_2)\big\|_{\mathbb{X}(t)}, \end{equation}
and the proof of the existence part of Theorem \ref{th:main1}  is achieved.  Moreover, the uniqueness part of Theorem \ref{th:main1} in $B(0,R)$ naturally follows.
\section{Time Decay Estimates}
\ \ \ \ \ In this section, we will establish the optimal time decay rates of the  global  strong solutions constructed in Theorem \ref{th:main1}.  We divide the proof into two steps.

\subsubsection*{Step 1: Low frequencies}

Denoting by $A(D)$ the semi-group associated to the system \eqref{equ:CTFS2} with $H_1\equiv H_2\equiv H_3\equiv H_4\equiv0$, for $U=((\sqrt{\beta_1}+\Lambda)c^{+},u^{+},(\sqrt{\beta_4}+\Lambda)c^{-},u^{-})$,   from  \eqref{gama8} and using Parseval's equality and the definition of $\ddq$,   we get
\begin{align*}\|e^{tA(D)}\ddq U\|_{L^2}&\lesssim  e^{-c_{0}2^{2q}t}\|\ddq U\|_{L^2}.
\end{align*}
Hence,   multiplying by $t^{\frac{N}{4}+\frac{s}{2}}2^{qs}$ and summing up on $q\leq q_{0}$, we readily have
\begin{equation}
\begin{split}
\label{low.2}
t^{\frac{N}{4}+\frac{s}{2}}\sum_{q\leq q_0}2^{qs}\|e^{tA(D)}\ddq U\|_{L^2}
&\lesssim
\sum_{q\leq q_0}2^{qs}e^{-c_{0}2^{2q}t}\|\ddq U\|_{L^2}t^{\frac{N}{4}+\frac{s}{2}}\\
&\lesssim
\sum_{q\leq q_0}2^{q(s+\frac{N}{2})}e^{-c_{0}2^{2q}t}\|\ddq U\|_{L^2}2^{q(-\frac{N}{2})}t^{\frac{N}{4}+\frac{s}{2}}\\
&\lesssim
\|U\|_{\dot B^{-\frac{N}{2}}_{2,\infty}}^\ell\sum_{q\leq q_0}2^{q(s+\frac{N}{2})}e^{-c_{0}2^{2q}t}t^{\frac{N}{4}+\frac{s}{2}}\\
&\lesssim
\|U\|_{\dot B^{-\frac{N}{2}}_{2,\infty}}^\ell\sum_{q\leq q_0}2^{q(s+\frac{N}{2})}e^{-c_{0}2^{2q}t}
t^{\frac{1}{2}(s+\frac{N}{2})}.
\end{split}
\end{equation}
As for any $\sigma>0$ there  exists a constant $C_\sigma$ so that
\begin{equation}\label{low.3}
\sup_{t\geq0}\sum_{q\in\mathbb{Z}}t^{\frac\sigma2}2^{q\sigma}e^{-c_{0}2^{2q}t}\leq C_\sigma.
\end{equation}
We get from \eqref{low.2} and \eqref{low.3} that for $s>-N/2,$
$$
\sup_{t\geq0}\, t^{\frac N4+\frac s2}\|e^{tA(D)}U\|_{\dot B^s_{2,1}}^\ell
\lesssim\|U\|_{\dot B^{-\frac{N}{2}}_{2,\infty}}^\ell.
$$
Furthermore, it is  obvious that  for $s>-N/2,$
$$
\|e^{tA(D)}U\|_{\dot B^s_{2,1}}^\ell
\lesssim \|U\|_{\dot B^{-\frac{N}{2}}_{2,\infty}}^\ell\sum_{q\leq q_0}2^{q(s+\frac{N}{2})}\lesssim\|U\|_{\dot B^{-\frac{N}{2}}_{2,\infty}}^\ell .
$$
Hence, setting $\langle t\rangle\eqdefa\sqrt{1+t^{2}}$, we get
\begin{equation}
\label{U}
\sup_{t\geq0}\, \langle t\rangle^{\frac N4+\frac s2}\|e^{tA(D)}U\|_{\dot B^s_{2,1}}^\ell
\lesssim\|U\|_{\dot B^{-\frac{N}{2}}_{2,\infty}}^\ell.
\end{equation}
Thus, from \eqref{U} and Duhamel's formula, we have
\begin{equation}
\begin{split}
&\big\|\big((\sqrt{\beta_1}+\Lambda)c^{+},u^{+},(\sqrt{\beta_4}+\Lambda)c^{-},u^{-}\big)\big\|_{\dot B^s_{2,1}}^\ell\\
&\quad\lesssim \sup_{t\geq0}\, \langle t\rangle^{-(\frac N4+\frac s2)}\big\|\big((\sqrt{\beta_1}+\Lambda)c^{+}_0,u^{+}_0,(\sqrt{\beta_4}+\Lambda)c^{-}_0,u^{-}_0\big)\big\|_{\dot B^{-\frac{N}{2}}_{2,\infty}}^\ell\\&\qquad+\int_{0}^{t}\langle t-\tau\rangle^{-(\frac N4+\frac s2)}\big\|\big((\sqrt{\beta_1}+\Lambda)H_{1},H_{2},(\sqrt{\beta_4}+\Lambda)H_{3},H_{4}\big)(\tau)
\big\|_{\dot B^{-\frac{N}{2}}_{2,\infty}}^\ell d\tau.
\end{split}
\end{equation}
We claim that for all $s\in(-N/2,2]$ and $t\geq0$,  then
\begin{equation}\label{s1234low1}\begin{split}
&\int_{0}^{t}\langle t-\tau\rangle^{-(\frac N4+\frac s2)}\big\|\big((\sqrt{\beta_1}+\Lambda)H_{1},H_{2},(\sqrt{\beta_4}+\Lambda)H_{3},H_{4}\big)(\tau)
\big\|_{\dot B^{-\frac{N}{2}}_{2,\infty}}^\ell d\tau
\\&\quad\lesssim\langle t\rangle^{-(\frac N4+\frac s2)}\big(1+X^{2}(t)\big)^{[\frac{N}{2}]+1}
\big(X^2(t)+D^2(t)+D^3(t)+D^4(t)\big),
\end{split}
\end{equation}
where $X(t)$ and $D(t)$ have been defined in \eqref{1.6} and \eqref{1.9}, respectively.\medbreak

Owing to the embedding $L^{1}\hookrightarrow\dot B^{-\frac{N}{2}}_{2,\infty}$,  it suffices to prove \eqref{s1234low1}  with $\|(H_{1},H_{2},H_{3},H_{4})(\tau)\|_{L^1}^{\ell}$ instead of $\|(H_{1},H_{2},H_{3},H_{4})(\tau)\|_{\dot  B^{-\frac{N}{2}}_{2,\infty}}^{\ell}$.

In  order to prove our claim,  we first present the following two important inequalities which will be frequently used in our process later.  By the smoothing effects of $(c^{+},u^{+},c^{-}, u^{-})$ in all frequencies,
 we have
\begin{equation}
\begin{split}
\label{6.6}
&\big\|\langle\tau\rangle^{\alpha}\big((\sqrt{\beta_1}+\Lambda)c^{+},u^{+},(\sqrt{\beta_4}+\Lambda)c^{-},u^{-}\big)\big\|_{\wt L^\infty_t(\dot B^{\frac N2-1}_{2,1})}^h
\\&\quad\lesssim\big\|\big((\sqrt{\beta_1}+\Lambda)c^{+},u^{+},(\sqrt{\beta_4}+\Lambda)c^{-},u^{-}\big)\big\|_{\wt L^\infty_t(\dot B^{\frac{N}{2}-1}_{2,1})}^h.
\\&\qquad+\big\|\tau^{\alpha}\big((\sqrt{\beta_1}+\Lambda)c^{+},u^{+},(\sqrt{\beta_4}+\Lambda)c^{-},u^{-}\big)\big\|_{\wt L^\infty_t(\dot B^{\frac N2+1}_{2,1})}^h
\\&\quad\lesssim X(t)+D(t).
\end{split}
\end{equation}
Further,  employing  a low-high decomposition and  \eqref{6.6}, for $N\geq2$ and $\alpha\geq\frac{N}{4}$, we obtain
\begin{equation}
\begin{split}
\label{6.7}&\sup_{0\leq \tau\leq t}\, \langle \tau\rangle^{\frac N4}\big\|\big((\sqrt{\beta_1}+\Lambda)c^{+},u^{+},(\sqrt{\beta_4}+\Lambda)c^{-},u^{-}\big)\big\|_{\dot B^{0}_{2,1}}
\\&\quad\lesssim \sup_{0\leq \tau\leq t}\, \langle \tau\rangle^{\frac N4}\big\|\big((\sqrt{\beta_1}+\Lambda)c^{+},u^{+},(\sqrt{\beta_4}+\Lambda)c^{-},u^{-}\big)\big\|_{\dot B^{0}_{2,1}}^\ell
\\&\qquad+\sup_{0\leq \tau\leq t}\, \langle \tau\rangle^{\frac N4}\big\|\big((\sqrt{\beta_1}+\Lambda)c^{+},u^{+},(\sqrt{\beta_4}+\Lambda)c^{-},u^{-}\big)\big\|_{\dot B^{0}_{2,1}}^{h}
\\&\quad\lesssim \sup_{0\leq \tau\leq t}\, \langle \tau\rangle^{\frac N4}\big\|\big((\sqrt{\beta_1}+\Lambda)c^{+},u^{+},(\sqrt{\beta_4}+\Lambda)c^{-},u^{-}\big)\big\|_{\dot B^{0}_{2,1}}^\ell
\\&\qquad+\sup_{0\leq \tau\leq t}\, \langle \tau\rangle^{\alpha}\big\|\big((\sqrt{\beta_1}+\Lambda)c^{+},u^{+},(\sqrt{\beta_4}+\Lambda)c^{-},u^{-}\big)\big\|_{\dot B^{\frac{N}{2}-1}_{2,1}}^{h}
\\&\quad\lesssim X(t)+D(t).
\end{split}
\end{equation}
To bound the term with $H_{1}$, we use the following decomposition:
$$ H_{1}=u^{+}\cdot\nabla c^{+}+ c^{+}\, \div (u^{+})^\ell + c^{+}\,\div (u^{+})^h.$$
Now, from H\"{o}lder's inequality,  the embedding $\dot B^{0}_{2,1}\hookrightarrow L^{2}$, the definitions of $D(t), \alpha$, Lemma \ref{lemma2.13}, \eqref{6.6} and \eqref{6.7}, one may write for all $s\in(\varepsilon-\frac N2,2]$,
\begin{equation}
\begin{split}\label{6.8}
&\int_0^t\langle t-\tau\rangle^{-(\frac N4+\frac s2)}\|(u^{+}\cdot\nabla c^{+})(\tau)\|_{L^1}\,d\tau\\
&\quad\lesssim\int_0^t\langle t-\tau\rangle^{-(\frac N4+\frac s2)}\|u^{+}\|_{L^2}\|\nabla c^{+}\|_{L^2}\,d\tau\\
&\quad\lesssim\int_0^t\langle t-\tau\rangle^{-(\frac N4+\frac s2)}\|u^{+}\|_{\dot  B^{0}_{2,1}}\|\nabla c^{+}\|_{\dot  B^{0}_{2,1}}\,d\tau\\
&\quad\lesssim\int_0^t\langle t-\tau\rangle^{-(\frac N4+\frac s2)}\|u^{+}\|_{\dot  B^{0}_{2,1}}\big(\|\nabla c^{+}\|_{\dot  B^{0}_{2,1}}^{\ell}+\|\nabla c^{+}\|_{\dot  B^{0}_{2,1}}^{h}\big)\,d\tau\\
&\quad\lesssim\big(\sup_{0\leq\tau\leq t}\langle \tau\rangle^{\frac N4}\|u^{+}(\tau)\|_{\dot  B^{0}_{2,1}}\big)
\big(\sup_{0\leq\tau\leq t}\langle \tau\rangle^{\frac N4+\frac 12}\|\nabla c^{+}(\tau)\|_{\dot  B^{0}_{2,1}}^{\ell}\big)
\int_0^t\langle t-\tau\rangle^{-(\frac N4+\frac s2)}\langle \tau\rangle^{-(\frac N2+\frac 12)}\,d\tau\\
&\qquad+\big(\sup_{0\leq\tau\leq t}\langle \tau\rangle^{\frac N4}\|u^{+}(\tau)\|_{\dot  B^{0}_{2,1}}\big)
\big(\sup_{0\leq\tau\leq t}\langle \tau\rangle^\alpha\|\nabla c^{+}(\tau)\|_{\dot  B^{\frac{N}{2}-1}_{2,1}}^h\big)
\int_0^t\langle t-\tau\rangle^{-(\frac N4+\frac s2)}\langle \tau\rangle^{-(\frac N4+\alpha)}\,d\tau\\
&\quad\lesssim \big(D^{2}(t)+X^{2}(t)\big)\int_0^t\langle t-\tau\rangle^{-(\frac N4+\frac s2)}\langle \tau\rangle^{-\min({\frac N2+\frac 12,\alpha+\frac{N}{4}})}\,d\tau\\
&\quad\lesssim \langle t\rangle^{-(\frac N4+\frac s2)} \big(D^{2}(t)+X^{2}(t)\big).
\end{split}
\end{equation}
The term  $c^{+}\,\div (u^{+})^\ell$ may be treated along the same lines,  and we have
\begin{equation}
\begin{split}
\label{6.9}
&\int_0^t\langle t-\tau\rangle^{-(\frac N4+\frac s2)}\|c^{+}\div (u^{+})^{\ell} \|_{L^1}\,d\tau\\
&\quad\lesssim\int_0^t\langle t-\tau\rangle^{-(\frac N4+\frac s2)}\|c^{+}\|_{L^2}\|\div (u^{+})^{\ell}\|_{L^2}\,d\tau\\
&\quad\lesssim\int_0^t\langle t-\tau\rangle^{-(\frac N4+\frac s2)}\|c^{+}\|_{\dot  B^{0}_{2,1}}\|\nabla u^{+}\|_{\dot  B^{0}_{2,1}}^{\ell}\,d\tau\\
&\quad\lesssim\big(\sup_{0\leq\tau\leq t}\langle \tau\rangle^{\frac N4}\|c^{+}(\tau)\|_{\dot  B^{0}_{2,1}}^{\ell}\big)
\big(\sup_{0\leq\tau\leq t}\langle \tau\rangle^{\frac{N}{4}+\frac 12}\| u^{+}(\tau)\|_{\dot  B^{1}_{2,1}}^\ell\big)
\int_0^t\langle t-\tau\rangle^{-(\frac N4+\frac s2)}\langle \tau\rangle^{-(\frac{N}{2}+\frac 12)}\,d\tau\\
&\quad\lesssim \langle t\rangle^{-(\frac N4+\frac s2)} \big(D^{2}(t)+X^{2}(t)\big).
\end{split}
\end{equation}
Regarding the term with $c^{+}\,\div (u^{+})^h,$ we get for all $t\geq2$ that
\begin{equation*}
\begin{split}
&\int_0^t\langle t-\tau\rangle^{-(\frac N4+\frac s2)}\|c^{+}\div (u^{+})^h(\tau)\|_{L^1}\,d\tau\\
&\quad\lesssim\int_0^t\langle t-\tau\rangle^{-(\frac N4+\frac s2)} \|c^{+}(\tau)\|_{\dot  B^{0}_{2,1}}\|\div u^{+}(\tau)\|^h_{\dot  B^{0}_{2,1}}\,d\tau\\
&\quad\lesssim\int_0^1\langle t-\tau\rangle^{-(\frac N4+\frac s2)} \|c^{+}(\tau)\|_{\dot  B^{0}_{2,1}}\|\div u^{+}(\tau)\|^h_{\dot  B^{0}_{2,1}}\,d\tau\\
&\qquad+\int_1^t\langle t-\tau\rangle^{-(\frac N4+\frac s2)} \|c^{+}(\tau)\|_{\dot  B^{0}_{2,1}}\|\div u^{+}(\tau)\|^h_{\dot  B^{0}_{2,1}}\,d\tau\\
&\quad\eqdefa I_{1}+I_{2}.
\end{split}
\end{equation*}
From the definitions of $X(t)$ and $D(t)$, we  obtain
\begin{align*}
I_{1}&\lesssim\langle t\rangle^{-(\frac N4+\frac s2)}
\sup _{0\leq\tau\leq1}\|c^{+}  (\tau)\|_{\dot  B^{0}_{2,1}}\int_0^1\|\div u^{+}(\tau)\|^h_{\dot  B^{0}_{2,1}}\,d\tau\\
&\lesssim\langle t\rangle^{-(\frac N4+\frac s2)}
\sup _{0\leq\tau\leq1}\|c^{+}(\tau)\|_{\dot  B^{0}_{2,1}}\int_0^1\| u^{+}(\tau)\|^h_{\dot  B^{\frac{N}{2}+1}_{2,1}}\,d\tau\\
&\lesssim\langle t\rangle^{-(\frac N4+\frac s2)} D(1)X(1),
\end{align*}
and,  using \eqref{6.6}, \eqref{6.7} and the fact that $\langle \tau\rangle\approx\tau$ when $\tau\geq1$, we get
\begin{align*}
I_{2} &\lesssim\int_1^t\langle t-\tau\rangle^{-(\frac N4+\frac s2)}
\|c^{+}(\tau)\|_{\dot  B^{0}_{2,1}}\|\div u^{+}(\tau)\|^h_{\dot  B^{0}_{2,1}}\,d\tau\\
&\lesssim\big(\sup_{1\leq\tau\leq t}\langle \tau\rangle^{\frac N4}\|c^{+}(\tau)\|_{\dot  B^{0}_{2,1}}\big)
\big(\sup_{1\leq\tau\leq t}\| \tau\nabla u^{+}(\tau)\|^{h}_{\dot  B^{\frac{N}{2}}_{2,1}}\big)
\int_1^t\langle t-\tau\rangle^{-(\frac N4+\frac s2)}\langle \tau\rangle^{-(\frac N4+1)}\,d\tau\\
&\lesssim \langle t\rangle^{-(\frac N4+\frac s2)} \big(D^{2}(t)+X^{2}(t)\big).
\end{align*}
Thus, for $t\geq2,$ we conclude that
\begin{align}
\label{6.10}
&\int_0^t\langle t-\tau\rangle^{-(\frac N4+\frac s2)}\|c^{+}\div (u^{+})^h(\tau)\|_{L^1}\,d\tau\nonumber\\
&\quad\lesssim\langle t\rangle^{-(\frac N4+\frac s2)} \big(D^{2}(t)+X^{2}(t)\big).
\end{align}
The case $t\leq2$ is obvious as $\langle t\rangle\approx1$ and
$\langle t-\tau\rangle\approx1$ for $0\leq\tau\leq t\leq 2$, and
\begin{equation}
\begin{split}
\label{6.11}
&\int_0^t\|c^{+}\,\div (u^{+})^h\|_{L^1}\,d\tau\\
&\quad\lesssim \|c^{+}\|_{L^\infty_t(L^2)}\|\div(u^{+})\|^h_{L_t^1(L^2)}\\
&\quad\lesssim \|c^{+}\|_{L^\infty_t(\dot  B^{0}_{2,1})}
\|\div (u^{+})\|^h_{L_t^1(\dot  B^{0}_{2,1})}\\
&\quad\lesssim \|c^{+}\|_{L^\infty_t(\dot  B^{0}_{2,1})}
\| u^{+}\|^h_{L_t^1(\dot  B^{\frac{N}{2}+1}_{2,1})}\\
&\quad\lesssim X(t)D(t).
\end{split}
\end{equation}
From \eqref{6.8}-\eqref{6.11}, we get
\begin{equation*}\label{s1}
\int_0^t\langle t-\tau\rangle^{-(\frac{N}{4}+\frac{s}{2})} \big\|H_{1}(\tau)\big\|_{\dot B^{-\frac{N}{2}}_{2,\infty}}^\ell d\tau
\lesssim\langle t\rangle^{-(\frac N4+\frac s2)}
\big(X^2(t)+D^2(t)\big).
\end{equation*}
The term  $H_{3}$ may be treated along the same lines,  and we obtain
\begin{equation*}\label{s3}
\int_0^t\langle t-\tau\rangle^{-(\frac{N}{4}+\frac{s}{2})} \big\|H_{3}(\tau)\big\|_{\dot B^{-\frac{N}{2}}_{2,\infty}}^\ell d\tau
\lesssim\langle t\rangle^{-(\frac N4+\frac s2)}
\big(X^2(t)+D^2(t)\big).
\end{equation*}
Next, we  bound  $H_{2}^{i}$ as follows. For the first part of $H_{2}^{i}$, employing \eqref{6.7} and Lemma \ref{p27} (i),  we write that
\begin{equation}
\begin{split}\label{6.12}
&\int_0^t\langle t-\tau\rangle^{-(\frac N4+\frac s2)}\|g_{+}(c^{+},c^{-})\partial_{i}c^{+}(\tau)\|_{L^1}\,d\tau\\
&\quad\lesssim\int_0^t\langle t-\tau\rangle^{-(\frac N4+\frac s2)}
\|g_{+}(c^{+},c^{-})\|_{L^2}\|\nabla c^{+}\|_{L^2}\,d\tau\\
&\quad\lesssim\int_0^t\langle t-\tau\rangle^{-(\frac N4+\frac s2)}\|g_{+}(c^{+},c^{-})\|_{\dot  B^{0}_{2,1}}\|\nabla c^{+}\|_{\dot  B^{0}_{2,1}}\,d\tau\\
&\quad\lesssim\big(1+X^{2}(t)\big)^{[\frac{N}{2}]+1}\int_0^t\langle t-\tau\rangle^{-(\frac N4+\frac s2)}\|(c^{+},c^{-})\|_{\dot  B^{0}_{2,1}}\|\nabla c^{+}\|_{\dot  B^{0}_{2,1}}\,d\tau\\
&\quad\lesssim\big(1+X^{2}(t)\big)^{[\frac{N}{2}]+1}\int_0^t\langle t-\tau\rangle^{-(\frac N4+\frac s2)}
\|(c^{+},c^{-})\|_{\dot B^{0}_{2,1}}
\big(\|\nabla c^{+}\|_{\dot  B^{0}_{2,1}}^{\ell}+\|\nabla c^{+}\|_{\dot  B^{0}_{2,1}}^{h}\big)\,d\tau\\
&\quad\lesssim\big(1+X^{2}(t)\big)^{[\frac{N}{2}]+1}\big(\sup_{0\leq\tau\leq t}\langle \tau\rangle^{\frac N4}\|(c^{+},c^{-})(\tau)\|_{\dot  B^{0}_{2,1}}\big)
\big(\sup_{0\leq\tau\leq t}\langle \tau\rangle^{\frac N4+\frac 12}\|\nabla c^{+}(\tau)\|_{\dot  B^{0}_{2,1}}^{\ell}\big)
\\&\qquad\times\int_0^t\langle t-\tau\rangle^{-(\frac N4+\frac s2)}\langle \tau\rangle^{-(\frac N2+\frac 12)}\,d\tau\\
&\qquad+\big(1+X^{2}(t)\big)^{[\frac{N}{2}]+1}\big(\sup_{0\leq\tau\leq t}\langle \tau\rangle^{\frac N4}\|(c^{+},c^{-})(\tau)\|_{\dot  B^{0}_{2,1}}\big)
\big(\sup_{0\leq\tau\leq t}\langle \tau\rangle^\alpha\|\nabla c^{+}(\tau)\|_{\dot  B^{\frac{N}{2}-1}_{2,1}}^h\big)
\\&\qquad\times\int_0^t\langle t-\tau\rangle^{-(\frac N4+\frac s2)}\langle \tau\rangle^{-(\alpha+\frac N4)}\,d\tau\\
&\quad\lesssim \big(1+X^{2}(t)\big)^{[\frac{N}{2}]+1}\big(D^{2}(t)+X^{2}(t)\big)\int_0^t\langle t-\tau\rangle^{-(\frac N4+\frac s2)}\langle \tau\rangle^{-\min({\frac N2+\frac 12,\alpha+\frac{N}{4}})}\,d\tau\\
&\quad\lesssim \langle t\rangle^{-(\frac N4+\frac s2)}\big(1+X^{2}(t)\big)^{[\frac{N}{2}]+1}\big(D^{2}(t)+X^{2}(t)\big),
\end{split}
\end{equation}
where $g_{+}$ stands for some smooth function vanishing at $0$.\\
Similarly,
\begin{equation}
\begin{split}\label{6.13}
&\int_0^t\langle t-\tau\rangle^{-(\frac N4+\frac s2)}\|\tilde{g}_{+}(c^{+},c^{-})\partial_{i}c^{-}\|_{L^1}\,d\tau\\
&\quad\lesssim \langle t\rangle^{-(\frac N4+\frac s2)} \big(1+X^{2}(t)\big)^{[\frac{N}{2}]+1}\big(D^{2}(t)+X^{2}(t)\big).
\end{split}
\end{equation}
To bound the term with $(u^{+}\cdot\nabla)u_{i}^{+}$, we employ the following decomposition:
$$(u^{+}\cdot\nabla)u_{i}^{+}=(u^{+}\cdot\nabla)(u_{i}^{+})^{\ell}
+(u^{+}\cdot\nabla)(u_{i}^{+})^{h}.
$$
By the same estimate of the term $c^{+}\,\div (u^{+})$, we have
\begin{align}\label{6.14}
&\int_0^t\langle t-\tau\rangle^{-(\frac N4+\frac s2)}
\|(u^{+}\cdot\nabla)(u_{i}^{+})(\tau)\|_{L^1}\,d\tau\nonumber\\
&\quad\lesssim\langle t\rangle^{-(\frac N4+\frac s2)} \big(D^{2}(t)+X^{2}(t)\big).
\end{align}
To deal with the term  $\mu^{+}h_{+}(c^{+},c^{-})\partial_{j}c^{+}\partial_{j}u^{+}_{i}$, we take the following decomposition:
$$\mu^{+}h_{+}(c^{+},c^{-})\partial_{j}c^{+}\partial_{j}u^{+}_{i}
=\mu^{+}h_{+}(c^{+},c^{-})\partial_{j}c^{+}\partial_{j}(u^{+}_{i})^{\ell}
+\mu^{+}h_{+}(c^{+},c^{-})\partial_{j}c^{+}\partial_{j}(u^{+}_{i})^{h}.
$$
For the term  $\mu^{+}h_{+}(c^{+},c^{-})\partial_{j}c^{+}\partial_{j}(u^{+}_{i})^{\ell}$, we have
\begin{equation*}
\begin{split}
\label{s24l}
&\int_0^t\langle t-\tau\rangle^{-(\frac N4+\frac s2)}
\|\mu^{+}h_{+}(c^{+},c^{-})\partial_{j}c^{+}\partial_{j}(u^{+}_{i})^{\ell}\|_{L^1}\,d\tau\\
&\quad\lesssim\int_0^t\langle t-\tau\rangle^{-(\frac N4+\frac s2)}
\|h_{+}(c^{+},c^{-})\|_{L^2}
\|\nabla c^{+}\nabla (u^{+})^{\ell}\|_{L^2}\,d\tau\\
&\quad\lesssim\int_0^t\langle t-\tau\rangle^{-(\frac N4+\frac s2)}
\|\nabla c^{+}\|_{\dot  B^{0}_{2,1}}\|\nabla(u^{+})^{\ell}\|_{\dot  B^{\frac N2}_{2,1}}\,d\tau\\
&\qquad+\int_0^t\langle t-\tau\rangle^{-(\frac N4+\frac s2)}\big(1+X^{2}(t)\big)^{[\frac{N}{2}]+1}
\|(c^{+},c^{-})\|_{\dot  B^{0}_{2,1}}
\|\nabla c^{+}\|_{\dot  B^{0}_{2,1}}\|\nabla(u^{+})^{\ell}\|_{\dot  B^{\frac N2}_{2,1}}\,d\tau\\
&\quad\eqdefa L_{1}+L_{2}.
\end{split}
\end{equation*}
We bound the two terms $L_{1}$ and $L_{2}$ as follows respectively, from \eqref{6.7}, we have
\begin{equation*}
\begin{split}
&L_{1}\lesssim\int_0^t\langle t-\tau\rangle^{-(\frac N4+\frac s2)}
\|\nabla c^{+}\|_{\dot  B^{0}_{2,1}}
\|\nabla u^{+}\|_{\dot  B^{0}_{2,1}}^{\ell}\,d\tau\\
&\quad\lesssim
\big(\sup_{0\leq\tau\leq t}\langle \tau\rangle^{\frac N4}\|\nabla c^{+}(\tau)\|_{\dot  B^{0}_{2,1}}\big)
\big(\sup_{0\leq\tau\leq t}\langle \tau\rangle^{\frac{N}{2}+\frac 12}\| \nabla u^{+}(\tau)\|_{\dot  B^{\frac N2}_{2,1}}^\ell\big)
\\&\qquad\times\int_0^t\langle t-\tau\rangle^{-(\frac N4+\frac s2)}\langle \tau\rangle^{-(\frac{3N}{4}+\frac{1}{2})}\,d\tau\\
&\quad\lesssim \langle t\rangle^{-(\frac N4+\frac s2)} \big(D^{2}(t)+X^{2}(t)\big),
\end{split}
\end{equation*}
and
\begin{equation*}
\begin{split}
&L_{2}\lesssim\int_0^t\langle t-\tau\rangle^{-(\frac N4+\frac s2)}\big(1+X^{2}(t)\big)^{[\frac{N}{2}]+1}
\|(c^{+},c^{-})\|_{\dot  B^{0}_{2,1}}
\|\nabla c^{+}\|_{\dot  B^{0}_{2,1}}
\|\nabla u^{+}\|_{\dot  B^{0}_{2,1}}^{\ell}\,d\tau\\
&\quad\lesssim\big(1+X^{2}(t)\big)^{[\frac{N}{2}]+1}\big(\sup_{0\leq\tau\leq t}\langle \tau\rangle^{\frac N4}\|(c^{+},c^{-})(\tau)\|_{\dot  B^{0}_{2,1}}\big)\langle \tau\rangle^{\frac N4}\| \nabla c^{+}(\tau)\|_{\dot  B^{0}_{2,1}}\langle \tau\rangle^{\frac{N}{2}+\frac 12}\| \nabla u^{+}(\tau)\|_{\dot  B^{\frac N2}_{2,1}}^\ell\big)
\\&\qquad\times\int_0^t\langle t-\tau\rangle^{-(\frac N4+\frac s2)}\langle \tau\rangle^{-(N+\frac 12)}\,d\tau\\
&\quad\lesssim \langle t\rangle^{-(\frac N4+\frac s2)} \big(1+X^{2}(t)\big)^{[\frac{N}{2}]+1}\big(D^{2}(t)+ D^{3}(t)+X^{4}(t)\big).
\end{split}
\end{equation*}
Thus,
\begin{equation}
\begin{split}
\label{s24l}
&\int_0^t\langle t-\tau\rangle^{-(\frac N4+\frac s2)}
\|\mu^{+}h_{+}(c^{+},c^{-})\partial_{j}c^{+}\partial_{j}(u^{+}_{i})^{\ell}\|_{L^1}\,d\tau\\
&\quad\lesssim\langle t\rangle^{-(\frac N4+\frac s2)} \big(1+X^{2}(t)\big)^{[\frac{N}{2}]+1}\big(D^{2}(t)+X^{2}(t)+ D^{3}(t)+X^{4}(t)\big).
\end{split}
\end{equation}
Regarding the term  $\mu^{+}h_{+}(c^{+},c^{-})\partial_{j}c^{+}\partial_{j}(u^{+}_{i})^{h},$ we also get,
\begin{equation*}
\begin{split}
&\int_0^t\langle t-\tau\rangle^{-(\frac N4+\frac s2)}
\|\mu^{+}h_{+}(c^{+},c^{-})\partial_{j}c^{+}\partial_{j}(u^{+}_{i})^{h}
(\tau)\|_{L^1}\,d\tau\\
&\quad\lesssim\int_0^t\langle t-\tau\rangle^{-(\frac N4+\frac s2)} \|h_{+}(c^{+},c^{-})(\tau)\|_{L^2}
\|\partial_{j}c^{+}\partial_{j}(u^{+}_{i})^{h}(\tau)\|_{L^2}\,d\tau\\
&\quad\lesssim\int_0^t\langle t-\tau\rangle^{-(\frac N4+\frac s2)} \big(1+\big(1+X^{2}(t)\big)^{[\frac{N}{2}]+1}\|(c^{+},c^{-})(\tau)\|_{\dot  B^{0}_{2,1}}\big)
\|\nabla c^{+}(\tau)\|_{\dot  B^{0}_{2,1}}
\|\nabla u^{+}(\tau)\|^h_{\dot  B^{\frac N2}_{2,1}}\,d\tau\\
&\quad\lesssim\int_0^t\langle t-\tau\rangle^{-(\frac N4+\frac s2)}
\|\nabla c^{+}(\tau)\|_{\dot  B^{0}_{2,1}}
\|\nabla u^{+}(\tau)\|^h_{\dot  B^{\frac N2}_{2,1}}\,d\tau\\
&\qquad+
\int_0^t\langle t-\tau\rangle^{-(\frac N4+\frac s2)} \big(1+X^{2}(t)\big)^{[\frac{N}{2}]+1}\|(c^{+},c^{-})(\tau)\|_{\dot  B^{0}_{2,1}}
\|\nabla c^{+}(\tau)\|_{\dot  B^{0}_{2,1}}
\|\nabla u^{+}(\tau)\|^h_{\dot  B^{\frac N2}_{2,1}}\,d\tau\\
&\quad\eqdefa M_{1}+M_{2}.
\end{split}
\end{equation*}
We deal with the two terms $M_{1}$ and $M_{2}$ in the following,
if $t\geq2,$
\begin{equation*}
\begin{split}
&M_{1}\lesssim\int_0^1\langle t-\tau\rangle^{-(\frac N4+\frac s2)}
\|\nabla c^{+}(\tau)\|_{\dot  B^{0}_{2,1}}
\|\nabla u^{+}(\tau)\|^h_{\dot  B^{\frac N2}_{2,1}}\,d\tau\\
&\qquad+\int_1^t\langle t-\tau\rangle^{-(\frac N4+\frac s2)}
\|\nabla c^{+}(\tau)\|_{\dot  B^{0}_{2,1}}
\|\nabla u^{+}(\tau)\|^h_{\dot  B^{\frac N2}_{2,1}}\,d\tau\\
&\quad\eqdefa M_{11}+M_{12}.
\end{split}
\end{equation*}
Using the definitions of $X(t)$ and $D(t)$, we  obtain
\begin{align*}
M_{11}\lesssim\langle t\rangle^{-(\frac N4+\frac s2)}
\sup _{0\leq\tau\leq1}\|\nabla c^{+}(\tau)\|_{\dot  B^{0}_{2,1}}
\int_0^1\|\nabla u^{+}(\tau)\|^h_{\dot B^{\frac{N}{2}}_{2,1}}\,d\tau\lesssim\langle t\rangle^{-(\frac N4+\frac s2)} D(1)X(1),
\end{align*}
and,  employing \eqref{6.7} and the fact that $\langle \tau\rangle\approx\tau$ when $\tau\geq1$, we have
\begin{align*}
M_{12}&\lesssim
\big(\sup_{1\leq\tau\leq t}\langle \tau\rangle^{\frac N4}
\|\nabla c^{+}(\tau)\|_{\dot  B^{0}_{2,1}}\big)
\big(\sup_{1\leq\tau\leq t}\| \tau\nabla u^{+}(\tau)\|^{h}_{\dot  B^{\frac{N}{2}}_{2,1}}\big)
\\&\qquad\times\int_1^t\langle t-\tau\rangle^{-(\frac N4+\frac s2)}\langle \tau\rangle^{-(\frac N4+1)}\,d\tau\\
&\lesssim \langle t\rangle^{-(\frac N4+\frac s2)} D^{2}(t).
\end{align*}
For the term  $M_{2}$, we have
\begin{equation*}
\begin{split}
M_{2}&\lesssim\int_0^1\langle t-\tau\rangle^{-(\frac N4+\frac s2)} \big(1+X^{2}(t)\big)^{[\frac{N}{2}]+1} \|(c^{+},c^{-})(\tau)\|_{\dot  B^{0}_{2,1}}
\|\nabla c^{+}(\tau)\|_{\dot  B^{0}_{2,1}}
\|\nabla u^{+}(\tau)\|^h_{\dot  B^{\frac N2}_{2,1}}\,d\tau\\
&\qquad+\int_1^t\langle t-\tau\rangle^{-(\frac N4+\frac s2)}  \big(1+X^{2}(t)\big)^{[\frac{N}{2}]+1}\|(c^{+},c^{-})(\tau)\|_{\dot  B^{0}_{2,1}}
\|\nabla c^{+}(\tau)\|_{\dot  B^{0}_{2,1}}
\|\nabla u^{+}(\tau)\|^h_{\dot  B^{\frac N2}_{2,1}}\,d\tau\\
&\eqdefa M_{21}+M_{22}.
\end{split}
\end{equation*}
Remembering the definitions of $X(t)$ and $D(t)$, we  obtain
\begin{align*}
M_{21}&\lesssim \big(1+X^{2}(t)\big)^{[\frac{N}{2}]+1}\langle t\rangle^{-(\frac N4+\frac s2)}
\sup _{0\leq\tau\leq1}\|(c^{+},c^{-})(\tau)\|_{\dot  B^{0}_{2,1}}
\sup _{0\leq\tau\leq1}\|\nabla c^{+}(\tau)\|_{\dot  B^{0}_{2,1}}
\int_0^1\|\nabla u^{+}(\tau)\|^h_{\dot B^{\frac{N}{2}}_{2,1}}\,d\tau\\
&\lesssim\langle t\rangle^{-(\frac N4+\frac s2)} \big(1+X^{2}(t)\big)^{[\frac{N}{2}]+1} \big(X^{2}(1)+D^{4}(1)\big),
\end{align*}
and
\begin{align*}
M_{22}&\lesssim\big(1+X^{2}(t)\big)^{[\frac{N}{2}]+1} \big(\sup_{1\leq\tau\leq t}\langle \tau\rangle^{\frac N4}
\|(c^{+},c^{-})(\tau)\|_{\dot  B^{0}_{2,1}}\sup_{1\leq\tau\leq t}\langle \tau\rangle^{\frac N4}
\|\nabla c^{+}(\tau)\|_{\dot  B^{0}_{2,1}} \sup_{1\leq\tau\leq t}\langle \tau\rangle^{\alpha}  \|\nabla u^{+}(\tau)\|^h_{\dot  B^{\frac N2}_{2,1}}  \big)
\\&\quad\times\int_1^t\langle t-\tau\rangle^{-(\frac N4+\frac s2)}
\langle \tau\rangle^{-(\frac N2+\alpha)}\,d\tau\\
&\lesssim \langle t\rangle^{-(\frac N4+\frac s2)}\big(1+X^{2}(t)\big)^{[\frac{N}{2}]+1} D^{3}(t).
\end{align*}
Therefore, for $t\geq2,$  we obatin
\begin{align}
\label{s24h2}
&\int_0^t\langle t-\tau\rangle^{-(\frac N4+\frac s2)}
\|\mu^{+}h_{+}(c^{+},c^{-})\partial_{j}c^{+}\partial_{j}(u^{+}_{i})^{h}(\tau)\|_{L^1}\,d\tau\nonumber\\
&\quad\lesssim\langle t\rangle^{-(\frac N4+\frac s2)}\big(1+X^{2}(t)\big)^{[\frac{N}{2}]+1} \big(X^{2}(t)+D^{2}(t)+D^{3}(t)+D^{4}(t)\big).
\end{align}
The case $t\leq2$ is obvious as $\langle t\rangle\approx1$ and
$\langle t-\tau\rangle\approx1$ for $0\leq\tau\leq t\leq 2$, and
\begin{equation}
\begin{split}
\label{s24h1}
&\int_0^t\|\mu^{+}h_{+}(c^{+},c^{-})\partial_{j}c^{+}\partial_{j}(u^{+}_{i})^{h}
(\tau)\|_{L^1}\,d\tau\\
&\quad\lesssim \|h_{+}(c^{+},c^{-})\|_{L^\infty_t(L^2)}
\|\nabla c^{+}\nabla (u^{+})^h\|_{L_t^1(L^2)}\\
&\quad\lesssim\big( \|h_{+}(c^{+},c^{-})-\frac{(\mathcal{C}^{2}\alpha^{-})(1,1)}
{s_{-}^{2}(1,1)}\|_{\dot  B^{0}_{2,1}}+1\big)
\|\nabla c^{+}\|_{L^\infty_t(\dot  B^{0}_{2,1})}
\|\nabla u^{+}\|^h_{L_t^1(\dot  B^{\frac{N}{2}}_{2,1})}\\
&\quad\lesssim \Big(1+\big(1+X^{2}(t)\big)^{[\frac{N}{2}]+1}\|(c^{+},c^{-})\|_{L^\infty_t(\dot  B^{0}_{2,1})}\Big)
\|\nabla c^{+}\|_{L^\infty_t(\dot  B^{0}_{2,1})}
\|u^{+}\|^h_{L_t^1(\dot  B^{\frac{N}{2}+1}_{2,1})}\\
&\quad\lesssim \big(1+X^{2}(t)\big)^{[\frac{N}{2}]+1}\big(X^{2}(t)+D^{2}(t)+D^{4}(t)\big).
\end{split}
\end{equation}
From \eqref{s24l}-\eqref{s24h1}, we finally conclude that
\begin{equation}
\begin{split}
\label{s24}
&\int_0^t\|\mu^{+}h_{+}(c^{+},c^{-})\partial_{j}c^{+}\partial_{j}u^{+}_{i}
(\tau)\|_{L^1}\,d\tau\\&\quad\lesssim \langle t\rangle^{-(\frac N4+\frac s2)}\big(1+X^{2}(t)\big)^{[\frac{N}{2}]+1} \big(X^{2}(t)+D^{2}(t)+D^{3}(t)+D^{4}(t)\big).
\end{split}
\end{equation}
Similarly,  we also obtain  the 	corresponding estimates of other terms
$\mu^{+}k_{+}(c^{+},c^{-})\partial_{j}c^{-}\partial_{j}u^{+}_{i},$\\
$\mu^{+}h_{+}(c^{+},c^{-})\partial_{j}c^{+}\partial_{i}u^{+}_{j},$
$\mu^{+}k_{+}(c^{+},c^{-})\partial_{j}c^{-}\partial_{i}u^{+}_{j},
\lambda^{+}h_{+}(c^{+},c^{-})\partial_{i}c^{+}\partial_{j}u^{+}_{j}~\text{and}~
\lambda^{+}k_{+}(c^{+},c^{-})\partial_{i}c^{-}\partial_{j}u^{+}_{j}$. Here, we omit the details.

From the low-high frequency  decomposition for $\mu^{+}l_{+}(c^{+},c^{-})\partial_{j}^{2}u_{i}^{+},$  we have
\begin{equation*}
\mu^{+}l_{+}(c^{+},c^{-})\partial_{j}^{2}u_{i}^{+}
=\mu^{+}l_{+}(c^{+},c^{-})\partial_{j}^{2}(u_{i}^{+})^\ell
+\mu^{+}l_{+}(c^{+},c^{-})\partial_{j}^{2}(u_{i}^{+})^h,
\end{equation*}
where $l_{+}$ stands for some smooth function vanishing at $0$.
Thus,
 \begin{equation}
 \begin{split}
 \label{s10l}
 &\int_0^t\langle t-\tau\rangle^{-(\frac N4+\frac s2)}\|\mu^{+}l_{+}(c^{+},c^{-})\partial_{j}^{2}(u_{i}^{+})^\ell\|_{L^1}\,d\tau\\
 &\quad\lesssim \big(1+X^{2}(t)\big)^{[\frac{N}{2}]+1} \big(\sup_{\tau\in[0,t]} \langle\tau\rangle^{\frac N4}\|(c^{+},c^{-})(\tau)\|_{\dot  B^{0}_{2,1}}\big)
 \big(\sup_{\tau\in[0,t]} \langle\tau\rangle^{\frac N4+1}
 \|\nabla^{2}u^{+}\|_{\dot  B^{0}_{2,1}}^\ell)\big)
\\&\qquad\times \int_0^t \langle t-\tau\rangle^{-(\frac N4+\frac s2)}
 \langle\tau\rangle^{-(\frac N2+1)}\,d\tau\\
 &\quad\lesssim \langle t\rangle^{-(\frac N4+\frac s2)}\big(1+X^{2}(t)\big)^{[\frac{N}{2}]+1}  D^{2}(t).
 \end{split}
 \end{equation}
To handle the term $\mu^{+}l_{+}(c^{+},c^{-})\partial_{j}^{2}(u_{i}^{+})^h$, we consider the cases $t\geq2$ and $t\leq2$ respectively. When $t\geq2$, then we have
 \begin{align*}
 &\int_0^t\langle t-\tau\rangle^{-(\frac N4+\frac s2)}\|\mu^{+}l_{+}(c^{+},c^{-})\partial_{j}^{2}(u_{i}^{+})^h\|_{L^1}\,d\tau
 \\&\quad\lesssim\int_0^t\langle t-\tau\rangle^{-(\frac N4+\frac s2)}\|l_{+}(c^{+},c^{-})\|_{\dot  B^{0}_{2,1}}
 \|\nabla^{2} u^{+}\|^h_{\dot  B^{0}_{2,1}}\,d\tau\\
 &\quad\lesssim\int_0^1\langle t-\tau\rangle^{-(\frac N4+\frac s2)}\big(1+X^{2}(t)\big)^{[\frac{N}{2}]+1} \|(c^{+},c^{-})\|_{\dot  B^{0}_{2,1}}
 \|\nabla^{2} u^{+}\|^h_{\dot  B^{0}_{2,1}}\,d\tau
\\&\qquad +\int_1^t\langle t-\tau\rangle^{-(\frac N4+\frac s2)}\big(1+X^{2}(t)\big)^{[\frac{N}{2}]+1} \|(c^{+},c^{-})\|_{\dot  B^{0}_{2,1}}\|\nabla^{2} u^{+}\|^h_{\dot  B^{0}_{2,1}}\,d\tau\\
 &\quad\eqdefa N_{1}+N_{2}.
 \end{align*}
From the definitions of $X(t)$ and $D(t)$, we  obtain
\begin{align*}
 N_{1}&=\int_0^1\langle t-\tau\rangle^{-(\frac N4+\frac s2)}\big(1+X^{2}(t)\big)^{[\frac{N}{2}]+1}
 \|(c^{+},c^{-})\|_{\dot  B^{0}_{2,1}}
 \|\nabla^{2} u^{+}\|^h_{\dot  B^{0}_{2,1}}\,d\tau\,d\tau\\
 &\lesssim\langle t\rangle^{-(\frac N4+\frac s2)}\big(1+X^{2}(t)\big)^{[\frac{N}{2}]+1}  \big(\sup_{\tau\in[0,1]}\|(c^{+},c^{-})\|_{\dot  B^{0}_{2,1}}\big)
\int_0^1 \|u^{+}\|^h_{\dot  B^{\frac{N}{2}+1}_{2,1}}\,d\tau\\
&\lesssim\langle t\rangle^{-(\frac N4+\frac s2)}\big(1+X^{2}(t)\big)^{[\frac{N}{2}]+1}  D(1)X(1),
\end{align*}
and, using the fact that $\langle \tau\rangle\approx\tau$ when $\tau\geq1$, from \eqref{6.7} we have
\begin{align*}
 N_{2}&=\int_1^t\langle t-\tau\rangle^{-(\frac N4+\frac s2)}
 \big(1+X^{2}(t)\big)^{[\frac{N}{2}]+1} \|(c^{+},c^{-})\|_{\dot  B^{0}_{2,1}}\|\nabla^{2} u^{+}\|^h_{\dot  B^{0}_{2,1}}\,d\tau\\
&\lesssim\big(1+X^{2}(t)\big)^{[\frac{N}{2}]+1} \big(\sup_{0\leq\tau\leq t}\langle \tau\rangle^{\frac N4}
\|(c^{+},c^{-})(\tau)\|_{\dot  B^{0}_{2,1}}\big)
\big(\sup_{0\leq\tau\leq t} \|\tau^{\alpha}\nabla u(\tau)\|_{\dot  B^{\frac{N}{2}}_{2,1}}^{h}\big)
\\&\quad\times\int_1^t\langle t-\tau\rangle^{-(\frac N4+\frac s2)}\langle \tau\rangle^{-(\frac N4+\alpha)}\,d\tau\\
&\lesssim \langle t\rangle^{-(\frac N4+\frac s2)}\big(1+X^{2}(t)\big)^{[\frac{N}{2}]+1}  \big(D^{2}(t)+X^{2}(t)\big).
\end{align*}
Thus, for $t\geq2,$ we arrive at
\begin{equation}
\label{s10h2}
\begin{split}
&\int_0^t\langle t-\tau\rangle^{-(\frac N4+\frac s2)}\|\mu^{+}l_{+}(c^{+},c^{-})\partial_{j}^{2}(u_{i}^{+})^h\|_{L^1}\,d\tau
\\&\quad\lesssim\langle t\rangle^{-(\frac N4+\frac s2)}\big(1+X^{2}(t)\big)^{[\frac{N}{2}]+1} (X^{2}(t)+D^{2}(t)).
\end{split}
\end{equation}
The case $t\leq2$ is obvious as $\langle t\rangle\approx1$ and $\langle t-\tau\rangle\approx1$ for $0\leq\tau\leq t\leq2,$
\begin{equation}
\label{s10h1}
\begin{split}
&\int_0^t\|\mu^{+}l_{+}(c^{+},c^{-})\partial_{j}^{2}(u_{i}^{+})^h\|_{L^1}\,d\tau\\
&\quad\lesssim\int_0^t\|l_{+}(c^{+},c^{-})\|_{\dot  B^{0}_{2,1}}
\|\nabla^{2}u^+\|_{\dot  B^{0}_{2,1}}^h\,d\tau\\
&\quad\lesssim\big(1+X^{2}(t)\big)^{[\frac{N}{2}]+1} \big(\sup_{\tau\in[0,1]}\|(c^{+},c^{-})(\tau)\|_{\dot  B^{0}_{2,1}}\big)
\int_0^1 \|u\|_{\dot  B^{\frac{N}{2}+1}_{2,1}}^h\,d\tau\\
&\quad\lesssim\big(1+X^{2}(t)\big)^{[\frac{N}{2}]+1}  D(t)X(t).
\end{split}
\end{equation}
From \eqref{s10l}-\eqref{s10h1}, we get
\begin{equation}
\label{s10}
\begin{split}
&\int_0^t\langle t-\tau\rangle^{-(\frac N4+\frac s2)}\|\mu^{+}l_{+}(c^{+},c^{-})\partial_{j}^{2}u_{i}^{+}\|_{L^1}\,d\tau
\\&\quad\lesssim\big(1+X^{2}(t)\big)^{[\frac{N}{2}]+1}  \langle t\rangle^{-(\frac N4+\frac s2)} \big(D^{2}(t)+X^{2}(t)\big).
\end{split}
\end{equation}
Similarly,
\begin{equation}
\label{s1last}
\begin{split}
&\int_0^t\langle t-\tau\rangle^{-(\frac N4+\frac s2)}\|(\mu^{+}+\lambda^{+})l_{+}(c^{+},c^{-})\partial_{i}\partial_{j}
u^{+}_{j}\|_{L^1}\,d\tau\\&\quad\lesssim \langle t\rangle^{-(\frac N4+\frac s2)}\big(1+X^{2}(t)\big)^{[\frac{N}{2}]+1}  \big(D^{2}(t)+X^{2}(t)\big).
\end{split}
\end{equation}
Thus,
\begin{equation*}\label{s2}\begin{split}
&\int_0^t\langle t-\tau\rangle^{-(\frac{N}{4}+\frac{s}{2})} \big\|H_{2}(\tau)\big\|_{\dot B^{-\frac{N}{2}}_{2,\infty}}^\ell d\tau
\\&\quad\lesssim\langle t\rangle^{-(\frac N4+\frac s2)}\big(1+X^{2}(t)\big)^{[\frac{N}{2}]+1}
\big(X^2(t)+D^2(t)+D^3(t)+D^4(t)\big).
\end{split}\end{equation*}
The term  $H_{4}$ may be treated along the same lines,  and we have
\begin{equation*}\label{s2}\begin{split}
&\int_0^t\langle t-\tau\rangle^{-(\frac{N}{4}+\frac{s}{2})} \big\|H_{4}(\tau)\big\|_{\dot B^{-\frac{N}{2}}_{2,\infty}}^\ell d\tau
\\&\quad\lesssim\langle t\rangle^{-(\frac N4+\frac s2)}\big(1+X^{2}(t)\big)^{[\frac{N}{2}]+1}
\big(X^2(t)+D^2(t)+D^3(t)+D^4(t)\big).
\end{split}\end{equation*}
Thus, we complete the proof of \eqref{s1234low1}.

Combining  with \eqref{U} and \eqref{s1234low1}, we conclude that for all $t\geq0$ and $s\in(\varepsilon-\frac{N}{2},2],$
\begin{equation}
\label{low}
\langle t\rangle^{\frac N4+\frac s2}
 \big\|\big(c^{+},\,u^{+},\,c^{-},\,u^{-}\big)\big\|^{\ell}_{\dot  B^{s}_{2,1}}\lesssim D_{0}+\big(1+X^{2}(t)\big)^{[\frac{N}{2}]+1}
\big(X^2(t)+D^2(t)+D^3(t)+D^4(t)\big).
\end{equation}
\subsubsection*{Step 2: High frequencies}
This step is devoted to bounding the last   term of $D(t)$. We first introduce the following system in terms of the weighted unknowns the term $(t^{\alpha}c^{+},t^{\alpha}u^{+},t^{\alpha}c^{-},t^{\alpha}u^{-} )$
\begin{equation}\label{equ:CTFS5}
\left\{
\begin{aligned}{}
&\p_t(t^{\alpha}c^{+})+\textrm{div}(t^{\alpha}u^{+})=\alpha t^{\alpha-1}c^{+}+  t^{\alpha}H_{1},\\
&\p_t(t^{\alpha}u^{+})+\beta_{1}\nabla (t^{\alpha}c^{+})
+\beta_{2}\nabla (t^{\alpha}c^{-})-\nu_{1}^{+}\Delta (t^{\alpha}u^{+})
-\nu_{2}^{+}\nabla\textrm{div}(t^{\alpha}u^{+})-\nabla \Delta (t^{\alpha}c^{+})\\&\qquad=\alpha t^{\alpha-1}u^{+}+t^{\alpha}H_{2},
\\&\p_t(t^{\alpha}c^{-})+\textrm{div}(t^{\alpha}u^{-})=\alpha t^{\alpha-1}c^{-}+t^{\alpha}H_{3},
\\&\p_t(t^{\alpha}u^{-})+\beta_{3}\nabla(t^{\alpha} c^{+})
+\beta_{4}\nabla(t^{\alpha} c^{-})-\nu_{1}^{-}\Delta(t^{\alpha} u^{-})
-\nu_{2}^{-}\nabla\textrm{div}(t^{\alpha}u^{-})-\nabla \Delta (t^{\alpha}c^{-})\\&\qquad=\alpha t^{\alpha-1}u^{-}+t^{\alpha}H_{4},\\
&(t^{\alpha}c^{+},t^{\alpha}u^{+},t^{\alpha}c^{-},t^{\alpha}u^{-})|_{t=0}=(0,0,0,0).
\end{aligned}
\right.
\end{equation}
From Lemma \ref{all part}, we have
\begin{equation}\label{regular estimate1}\begin{split}
&\big\|\tau^{\alpha}\big((\sqrt{\beta_1}+\Lambda)c^{+},u^{+},(\sqrt{\beta_4}+\Lambda)c^{-},u^{-}\big)\big\|
_{\wt L^\infty_t(\dot B^{\frac{N}{2}+1}_{2,1})}^{h}
\\&\quad\lesssim\big\|(\sqrt{\beta_1}+\Lambda)(\alpha\tau^{\alpha-1}c^{+}+\tau^{\alpha}H_{1})
\big\|_{\tilde{L}_{t}^{\infty}(\dot{B}^{\frac{N}{2}-1}_{2,1})}^{h}
\\&\qquad+\big\|\alpha\tau^{\alpha-1}u^{+}+\tau^{\alpha}H_{2}
\big\|_{\tilde{L}_{t}^{\infty}(\dot{B}^{\frac{N}{2}-1}_{2,1})}^{h}
\\&\qquad+\big\|(\sqrt{\beta_4}+\Lambda)(\alpha\tau^{\alpha-1}c^{-}+\tau^{\alpha}H_{3})
\big\|_{\tilde{L}_{t}^{\infty}(\dot{B}^{\frac{N}{2}-1}_{2,1})}^{h}
\\&\qquad+\big\|\alpha\tau^{\alpha-1}u^{-}+\tau^{\alpha}H_{4}
\big\|_{\tilde{L}_{t}^{\infty}(\dot{B}^{\frac{N}{2}-1}_{2,1})}^{h}.
\end{split}\end{equation}
We now handle the lower order linear terms on the right hand-side of the above, for $v\in \{(\sqrt{\beta_1}+\Lambda)c^{+},u^{+},(\sqrt{\beta_4}+\Lambda)c^{-},u^{-}\}$,
we have $$\|\alpha\tau^{\alpha-1}v
\|_{\tilde{L}_{t}^{\infty}(\dot{B}^{\frac{N}{2}-1}_{2,1})}^{h}\lesssim \|v
\|_{\tilde{L}_{t}^{\infty}(\dot{B}^{\frac{N}{2}-1}_{2,1})}^{h}\lesssim X(t),\, \hbox{for} \quad 0\leq \tau \leq t\leq2.$$
 When $t\geq2$, for  $0\leq \tau \leq 1$, we get
 $$\|\alpha\tau^{\alpha-1}v
\|_{\tilde{L}^{\infty}([0,1];\dot{B}^{\frac{N}{2}-1}_{2,1})}^{h}\lesssim \|v
\|_{\tilde{L}_{t}^{\infty}(\dot{B}^{\frac{N}{2}-1}_{2,1})}^{h}\lesssim X(t).$$
 When $t\geq2$, for  $1\leq \tau \leq t$, we have
\begin{equation}\label{regular estimate11}\begin{split}
\|\alpha\tau^{\alpha-1}v
\|_{\tilde{L}^{\infty}([1,t];\dot{B}^{\frac{N}{2}-1}_{2,1})}^{h}&=\alpha\sum_{j\geq j_0}2^{j(\frac{N}{2}+1)}2^{-2j}\|\tau^{\alpha-1}\dot{\Delta}_{j} v\|_{L^{\infty}([1,t],L^{2})}
\\&\lesssim \alpha 2^{-2j_0}\sum_{j\geq j_0}2^{j(\frac{N}{2}+1)}\|\tau^{\alpha}\dot{\Delta}_{j} v\|_{L^{\infty}([1,t],L^{2})}
\\&\lesssim \alpha 2^{-2j_0}\|\tau^{\alpha}v
\|_{\tilde{L}^{\infty}([1,t];\dot{B}^{\frac{N}{2}+1}_{2,1})}^{h}.
\end{split}\end{equation}
Choosing $j_0$ large enough such that  $$C\alpha 2^{-2j_0}\leq\frac{1}{4},$$ which implies that \eqref{regular estimate11} may be absorbed by  the  left hand-side of \eqref{regular estimate1}.
Thus
\begin{equation}\label{regular estimate111}\begin{split}
&\big\|\tau^{\alpha}\big((\sqrt{\beta_1}+\Lambda)c^{+},u^{+},(\sqrt{\beta_4}+\Lambda)c^{-},u^{-}\big)\big\|
_{\wt L^\infty_t(\dot B^{\frac{N}{2}+1}_{2,1})}^{h}
\\&\quad\lesssim  X(t)+\big\|\tau^{\alpha}(H_{1},H_{2},H_{3},H_{4})
\big\|_{\tilde{L}_{t}^{\infty}(\dot{B}^{\frac{N}{2}-1}_{2,1})}^{h}.
\end{split}\end{equation}
It now comes down to estimating the above nonlinear terms.
We first show the following inequalities which are repeatedly used later. For $\alpha=\frac{1}{2}(N+1-\varepsilon)$, we have
\begin{equation}\label{regular estimate1111}\begin{split}
\|\tau^{\alpha}\nabla(c^{+},c^{-})\|
_{\wt L^\infty_t(\dot B^{\frac{N}{2}}_{2,1})}&\lesssim  \|\tau^{\alpha}\nabla(c^{+},c^{-})\|
_{\wt L^\infty_t(\dot B^{\frac{N}{2}}_{2,1})}^\ell+\|\tau^{\alpha}\nabla(c^{+},c^{-})\|
_{\wt L^\infty_t(\dot B^{\frac{N}{2}}_{2,1})}^{h}\\
&\lesssim  \|\tau^{\alpha}(c^{+},c^{-})\|
_{L^\infty_t(\dot B^{\frac{N}{2}+1-\varepsilon}_{2,1})}^\ell+\|\tau^{\alpha}\nabla(c^{+},c^{-})\|
_{\wt L^\infty_t(\dot B^{\frac{N}{2}}_{2,1})}^{h}
\\&\lesssim D(t),
\end{split}\end{equation}
\begin{equation}\label{regular estimate11111}\begin{split}
\|\tau^{\alpha}(u^{+},u^{-})\|
_{\wt L^\infty_t(\dot B^{\frac{N}{2}+1}_{2,1})}&\lesssim  \|\tau^{\alpha}(u^{+},u^{-})\|
_{\wt L^\infty_t(\dot B^{\frac{N}{2}+1}_{2,1})}^\ell+\|\tau^{\alpha}(u^{+},u^{-})\|
_{\wt L^\infty_t(\dot B^{\frac{N}{2}+1}_{2,1})}^{h}\\
&\lesssim  \|\tau^{\alpha}(u^{+},u^{-})\|
_{L^\infty_t(\dot B^{\frac{N}{2}+1-\varepsilon}_{2,1})}^\ell+\|\tau^{\alpha}\nabla(u^{+},u^{-})\|
_{\wt L^\infty_t(\dot B^{\frac{N}{2}+1}_{2,1})}^{h}
\\&\lesssim D(t).
\end{split}\end{equation}
For $\|\tau^{\alpha}H_{1}
\|_{\tilde{L}_{t}^{\infty}(\dot{B}^{\frac{N}{2}-1}_{2,1})}^{h}$, from Proposition \ref{p26}, \eqref{regular estimate1111} and \eqref{regular estimate11111}, we have
\begin{equation}\label{regular estimate6.54}\begin{split}
\|\tau^{\alpha}H_{1}
\|_{\tilde{L}_{t}^{\infty}(\dot{B}^{\frac{N}{2}-1}_{2,1})}^{h}&\lesssim  \|\tau^{\alpha}c^{+}\textrm{div}u^{+}\|
_{\wt L^\infty_t(\dot B^{\frac{N}{2}-1}_{2,1})}^h+\|\tau^{\alpha}u^{+}\cdot \nabla c^{+})\|
_{\wt L^\infty_t(\dot B^{\frac{N}{2}-1}_{2,1})}^{h}\\
&\lesssim \|\tau^{\alpha}c^{+}\textrm{div}u^{+}\|
_{\wt L^\infty_t(\dot B^{\frac{N}{2}}_{2,1})}^h+\|\tau^{\alpha}u^{+}\cdot \nabla c^{+})\|
_{\wt L^\infty_t(\dot B^{\frac{N}{2}-1}_{2,1})}^{h}\\
&\lesssim \|c^{+}\|_{\wt L^\infty_t(\dot B^{\frac{N}{2}}_{2,1})} \|\tau^{\alpha}u^{+}\|
_{\wt L^\infty_t(\dot B^{\frac{N}{2}+1}_{2,1})}^h +\|u^{+}\|
_{\wt L^\infty_t(\dot B^{\frac{N}{2}-1}_{2,1})}\|\tau^{\alpha}\nabla c^{+}\|
_{\wt L^\infty_t(\dot B^{\frac{N}{2}}_{2,1})}
\\&\lesssim X(t)D(t).
\end{split}\end{equation}
Similarly,
\begin{equation}\label{regular estimate6.55}\begin{split}
\|\tau^{\alpha}H_{3}
\|_{\tilde{L}_{t}^{\infty}(\dot{B}^{\frac{N}{2}-1}_{2,1})}^{h}\lesssim X(t)D(t).
\end{split}\end{equation}
Next, we bound the term $\|\tau^{\alpha}H_{2}
\|_{\tilde{L}_{t}^{\infty}(\dot{B}^{\frac{N}{2}-1}_{2,1})}^{h}$ as follows. To bound the first part of $H_{2}^{i}$, employing \eqref{regular estimate1111}, Proposition \ref{p26} and  Lemma \ref{p27} (i),  we write that
\begin{equation}
\begin{split}\label{regular estimate6.56}
&\|\tau^{\alpha}g_{+}(c^{+},c^{-})\partial_{i}c^{+}
\|_{\tilde{L}_{t}^{\infty}(\dot{B}^{\frac{N}{2}-1}_{2,1})}^{h}\\
&\quad\lesssim\|g_{+}(c^{+},c^{-})\|_{\tilde{L}_{t}^{\infty}(\dot{B}^{\frac{N}{2}-1}_{2,1})}\|\tau^{\alpha}\partial_{i}c^{+}
\|_{\tilde{L}_{t}^{\infty}(\dot{B}^{\frac{N}{2}}_{2,1})}\\
&\quad\lesssim \big(1+X^{2}(t)\big)^{[\frac{N}{2}]+1}\|(c^{+},c^{-})\|_{\tilde{L}_{t}^{\infty}(\dot{B}^{\frac{N}{2}-1}_{2,1})}
\|\tau^{\alpha}\nabla c^{+}
\|_{\tilde{L}_{t}^{\infty}(\dot{B}^{\frac{N}{2}}_{2,1})}\\
&\quad\lesssim \big(1+X^{2}(t)\big)^{[\frac{N}{2}]+1}\big(D^{2}(t)+X^{2}(t)\big).
\end{split}
\end{equation}
Similar to  \eqref{regular estimate6.56}, we have
\begin{equation}
\begin{split}\label{regular estimate6.57}
&\|\tau^{\alpha}\tilde{g}_{+}(c^{+},c^{-})\partial_{i}c^{-}
\|_{\tilde{L}_{t}^{\infty}(\dot{B}^{\frac{N}{2}-1}_{2,1})}^{h}\\
&\quad\lesssim \big(1+X^{2}(t)\big)^{[\frac{N}{2}]+1}\big(D^{2}(t)+X^{2}(t)\big).
\end{split}
\end{equation}
To bound the term with $(u^{+}\cdot\nabla)u_{i}^{+}$, from Proposition \ref{p26} and \eqref{regular estimate11111}, we get
\begin{equation}
\begin{split}\label{regular estimate6.58}
&\|\tau^{\alpha}(u^{+}\cdot\nabla)u_{i}^{+}
\|_{\tilde{L}_{t}^{\infty}(\dot{B}^{\frac{N}{2}-1}_{2,1})}^{h}\\
&\quad\lesssim \|u^{+}
\|_{\tilde{L}_{t}^{\infty}(\dot{B}^{\frac{N}{2}-1}_{2,1})}\|\tau^{\alpha}\nabla u^{+}
\|_{\tilde{L}_{t}^{\infty}(\dot{B}^{\frac{N}{2}}_{2,1})}\\
&\quad\lesssim D^{2}(t)+X^{2}(t).
\end{split}
\end{equation}
Using \eqref{regular estimate11111}, Proposition \ref{p26} and  Lemma \ref{p27} (i),  we deduce that
 \begin{equation}
\begin{split}\label{regular estimate6.59}
&\|\tau^{\alpha}\mu^{+}h_{+}(c^{+},c^{-})\partial_{j}c^{+}\partial_{j}u^{+}_{i}
\|_{\tilde{L}_{t}^{\infty}(\dot{B}^{\frac{N}{2}-1}_{2,1})}^{h}\\
&\quad\lesssim\big(1+X^{2}(t)\big)^{[\frac{N}{2}]+1}\|(c^{+},c^{-})\|_{\tilde{L}_{t}^{\infty}(\dot{B}^{\frac{N}{2}}_{2,1})} \|\nabla c^{+}
\|_{\tilde{L}_{t}^{\infty}(\dot{B}^{\frac{N}{2}-1}_{2,1})}\|\tau^{\alpha}\nabla u^{+}
\|_{\tilde{L}_{t}^{\infty}(\dot{B}^{\frac{N}{2}}_{2,1})}\\
&\quad\lesssim\big(1+X^{2}(t)\big)^{[\frac{N}{2}]+1}X(t)\big( D^{2}(t)+X^{2}(t)\big).
\end{split}
\end{equation}
Similarly,  we also obtain  the 	corresponding estimates of other terms
$\mu^{+}k_{+}(c^{+},c^{-})\partial_{j}c^{-}\partial_{j}u^{+}_{i},$\\
$\mu^{+}h_{+}(c^{+},c^{-})\partial_{j}c^{+}\partial_{i}u^{+}_{j},$
$\mu^{+}k_{+}(c^{+},c^{-})\partial_{j}c^{-}\partial_{i}u^{+}_{j},
\lambda^{+}h_{+}(c^{+},c^{-})\partial_{i}c^{+}\partial_{j}u^{+}_{j}~\text{and}~
\lambda^{+}k_{+}(c^{+},c^{-})\partial_{i}c^{-}\partial_{j}u^{+}_{j}$. Here, we omit the details.

From \eqref{regular estimate11111}, Proposition \ref{p26} and  Lemma \ref{p27} (i),  we have
\begin{equation}
\begin{split}\label{regular estimate6.60}
&\|\tau^{\alpha}\mu^{+}l_{+}(c^{+},c^{-})\partial_{j}^{2}u_{i}^{+}
\|_{\tilde{L}_{t}^{\infty}(\dot{B}^{\frac{N}{2}-1}_{2,1})}^{h}\\
&\quad\lesssim\big(1+X^{2}(t)\big)^{[\frac{N}{2}]+1}\|(c^{+},c^{-})\|_{\tilde{L}_{t}^{\infty}(\dot{B}^{\frac{N}{2}}_{2,1})} \|\tau^{\alpha}\nabla^{2} u^{+}
\|_{\tilde{L}_{t}^{\infty}(\dot{B}^{\frac{N}{2}-1}_{2,1})}\\
&\quad\lesssim\big(1+X^{2}(t)\big)^{[\frac{N}{2}]+1}\big( D^{2}(t)+X^{2}(t)\big).
\end{split}
\end{equation}
Similarly,
\begin{equation}
\begin{split}\label{regular estimate6.611}
&\|\tau^{\alpha}(\mu^{+}+\lambda^{+})l_{+}(c^{+},c^{-})\partial_{i}\partial_{j}
u^{+}_{j}
\|_{\tilde{L}_{t}^{\infty}(\dot{B}^{\frac{N}{2}-1}_{2,1})}^{h}\\
&\quad\lesssim\big(1+X^{2}(t)\big)^{[\frac{N}{2}]+1}\|(c^{+},c^{-})\|_{\tilde{L}_{t}^{\infty}(\dot{B}^{\frac{N}{2}}_{2,1})} \|\tau^{\alpha}\nabla^{2} u^{+}
\|_{\tilde{L}_{t}^{\infty}(\dot{B}^{\frac{N}{2}-1}_{2,1})}\\
&\quad\lesssim\big(1+X^{2}(t)\big)^{[\frac{N}{2}]+1}\big( D^{2}(t)+X^{2}(t)\big).
\end{split}
\end{equation}
Combining with  \eqref{regular estimate6.56} and \eqref{regular estimate6.611}, we have
Thus,
\begin{equation}
\begin{split}\label{regular estimate6.612}
&\|\tau^{\alpha}H_2
\|_{\tilde{L}_{t}^{\infty}(\dot{B}^{\frac{N}{2}-1}_{2,1})}^{h}\\
&\quad\lesssim\big(1+X^{2}(t)\big)^{[\frac{N}{2}]+1}
\big(X^2(t)+D^2(t)+D^3(t)+D^4(t)\big).
\end{split}
\end{equation}
The term  $H_{4}$ may be treated along the same lines,  and we have
\begin{equation}
\begin{split}\label{regular estimate6.613}
&\|\tau^{\alpha}H_4
\|_{\tilde{L}_{t}^{\infty}(\dot{B}^{\frac{N}{2}-1}_{2,1})}^{h}\\
&\quad\lesssim\big(1+X^{2}(t)\big)^{[\frac{N}{2}]+1}
\big(X^2(t)+D^2(t)+D^3(t)+D^4(t)\big).
\end{split}
\end{equation}
Adding up  \eqref{regular estimate6.54}, \eqref{regular estimate6.55},   \eqref{regular estimate6.612} and \eqref{regular estimate6.613} to \eqref{regular estimate111}  yields
\begin{equation}\label{regular estimate6.614}\begin{split}
&\big\|\tau^{\alpha}\big((\sqrt{\beta_1}+\Lambda)c^{+},u^{+},(\sqrt{\beta_4}+\Lambda)c^{-},u^{-}\big)\big\|
_{\wt L^\infty_t(\dot B^{\frac{N}{2}+1}_{2,1})}^{h}
\\&\quad\lesssim  X(t)+\big(1+X^{2}(t)\big)^{[\frac{N}{2}]+1}
\big(X^2(t)+D^2(t)+D^3(t)+D^4(t)\big).
\end{split}\end{equation}
Finally, combining with \eqref{low} and  \eqref{regular estimate6.614},  for all $t\geq0,$ we have
\begin{align*}D(t)\lesssim  D_{0}+X(t)+\big(1+X^{2}(t)\big)^{[\frac{N}{2}]+1}
\big(X^2(t)+D^2(t)+D^3(t)+D^4(t)\big).
\end{align*}
As Theorem \ref{th:main1} ensures that $X(t)\lesssim  X(0)$  with $X(0) $ being small, and $X(0)^{\ell}=\|( c_{0}^{+}, u_{0}^{+}, c_{0}^{-}, u_{0}^{-})\|^{\ell}_{\dot B^{\frac N2-1}_{2,1}}\lesssim \|( c_{0}^{+}, u_{0}^{+}, c_{0}^{-}, u_{0}^{-})\|^{\ell}_{\dot B^{-\frac N2}_{2,\infty}}$,   one can conclude that \eqref{1.8} is fulfilled for all time if $D_{0}$ and $X(0) $ are small enough. This completes the proof of Theorem \ref{th:decay}.

\section{More Decay Estimates}
\ \ \ \ \ In this section, we present  some corollaries of Theorem \ref{th:decay}, which implies that the standard optimal   $L^{q}$-$L^{r}$ time decay rates of $(R^{+}-1,\,u^{+},\,R^{-}-1,\,u^{-})$.

\begin{Corollary}\label{cor1.1}  The solution $(R^{+}-1,\,u^{+},\,R^{-}-1,\,u^{-})$  constructed in  Theorem \ref{th:main1} satisfies
$$\displaylines{
\big\|\Lambda^{s}\big(R^{+}-1,R^{-}-1\big)\big\|_{L^2}
\lesssim \Big(D_0+\big\|\big(\nabla R^{+}_{0},\,u^{+}_{0},\,\nabla R^{-}_{0},\,u^{-}_{0}\big)\big\|^h_{\dot B^{\frac N2-1}_{2,1}}\Big)\langle t\rangle^{-\frac N4-\frac s2}\cr
\quad  \hbox{ if } \ \ -N/2<s\leq \min\{2,N/2 \},\cr
\big\|\Lambda^{s}\big(u^{+},u^{-}\big)\big\|_{L^2}
\lesssim \Big(D_0+\big\|\big(\nabla R^{+}_{0},\,u^{+}_{0},\,\nabla R^{-}_{0},\,u^{-}_{0}\big)\big\|^h_{\dot B^{\frac N2-1}_{2,1}}\Big)\langle t\rangle^{-\frac N4-\frac s2}\cr
\quad \hbox{ if } \ \ -N/2<s\leq \min\{2,N/2-1 \},}
$$
 where the fractional derivative
 operator $\Lambda^{\ell}$ is defined by $\Lambda^{\ell}f\triangleq\mathcal{F}^{-1}(|\cdot|^{\ell}\mathcal{F}f)$.
\end{Corollary}

\noindent{\bf Proof.} For the solution $(R^{+}-1,\,u^{+},\,R^{-}-1,\,u^{-})$  constructed in  Theorem \ref{th:main1}, applying  to homogeneous Littlewood-Paley decomposition
 for $R^{+}-1$, we have
$$\|\Lambda^s (R^{+}-1)\|_{L^2}\lesssim  \sum_{q\in \mathbb{Z}} \|\Delta_q\Lambda^s (R^{+}-1)\|_{L^2}=\|\Lambda^s (R^{+}-1)\|_{\dot B^0_{2,1}}.$$
Based on  Bernstein's inequalities and the low-high frequencies decomposition, we may write
$$
\sup_{t\in[0,T]} \langle t\rangle^{\frac N4+\frac s2}\|\Lambda^s (R^{+}-1)\|_{\dot B^0_{2,1}}\lesssim
 \big\|\langle t\rangle^{\frac N4+\frac s2}(R^{+}-1)\big\|_{L^\infty_T(\dot B^s_{2,1})}^\ell
 +  \big\|\langle t\rangle^{\frac N4+\frac s2}(R^{+}-1)\big\|_{L^\infty_T(\dot B^s_{2,1})}^h.
$$
If follows from Inequality \eqref{1.8} and  definitions of $D(t)$ and $\alpha$ that
$$
  \big\|\langle t\rangle^{\frac N4+\frac s2}(R^{+}-1)\big\|_{L^\infty_T(\dot B^s_{2,1})}^\ell\lesssim D_0+\big\|\big(\nabla R^{+}_{0},\,u^{+}_{0},\,\nabla R^{-}_{0},\,u^{-}_{0}\big)\big\|^h_{\dot B^{\frac N2-1}_{2,1}}
\quad\hbox{if }\ \  -N/2<s\leq 2.
$$ On the other hand,
for  $\alpha\geq\frac N4+\frac s2$,  $s\leq \min\{2,N/2 \},$ by \eqref{6.6}  we have
\begin{equation*}
\begin{split}
  \big\|\langle t\rangle^{\frac N4+\frac s2}(R^{+}-1)\big\|_{L^\infty_T(\dot B^s_{2,1})}^h
  &\lesssim  \big\|\langle t\rangle^{\alpha}(R^{+}-1)\big\|_{\tilde{L}^\infty_T(\dot B^{\frac{N}{2}}_{2,1})}^h\\
   &\lesssim X(t)+D(t) \quad\hbox{if }\ \  s\leq \min\{2,N/2 \}.
  \end{split}
\end{equation*}
Furthermore, from Theorem \ref{th:main1} we have $X(t)\lesssim  X(0)$, and $X(0)^{\ell}=\big\|( c_{0}^{+}, u_{0}^{+}, c_{0}^{-}, u_{0}^{-})\big\|^{\ell}_{\dot B^{\frac N2-1}_{2,1}}\lesssim \big\|( c_{0}^{+}, u_{0}^{+}, c_{0}^{-}, u_{0}^{-})\big\|^{\ell}_{\dot B^{-\frac N2}_{2,\infty}}$.  Combining with \eqref{1.8}, we get
\begin{equation*}
\begin{split}
  \big\|\langle t\rangle^{\frac N4+\frac s2}(R^{+}-1)\big\|_{L^\infty_T(\dot B^s_{2,1})}^h
  \lesssim  D_0+\big\|\big(\nabla R^{+}_{0},\,u^{+}_{0},\,\nabla R^{-}_{0},\,u^{-}_{0}\big)\big\|^h_{\dot B^{\frac N2-1}_{2,1}}\quad\hbox{if }\ \  s\leq \min\{2,N/2 \}.
  \end{split}
\end{equation*}
Thus, we obtain  the following desired result for $R^{+}-1$
$$\big\|\Lambda^{s}(R^{+}-1)\big\|_{L^2}
\lesssim \Big(D_0+\big\|\big(\nabla R^{+}_{0},\,u^{+}_{0},\,\nabla R^{-}_{0},\,u^{-}_{0}\big)\big\|^h_{\dot B^{\frac N2-1}_{2,1}}\Big)\langle t\rangle^{-\frac N4-\frac s2}
\quad  \hbox{ if } \ \ -N/2<s\leq \min\{2,N/2 \}.$$
Similarly,
$$\displaylines{
\|\Lambda^{s}u^{+}\|_{L^2}
\lesssim \Big(D_0+\big\|\big(\nabla R^{+}_{0},\,u^{+}_{0},\,\nabla R^{-}_{0},\,u^{-}_{0}\big)\big\|^h_{\dot B^{\frac N2-1}_{2,1}}\Big)\langle t\rangle^{-\frac N4-\frac s2}
\quad \hbox{ if } \ \ -N/2<s\leq \min\{2,N/2-1 \},\cr
\|\Lambda^{s}(R^{-}-1)\|_{L^2}
\lesssim \Big(D_0+\big\|\big(\nabla R^{+}_{0},\,u^{+}_{0},\,\nabla R^{-}_{0},\,u^{-}_{0}\big)\big\|^h_{\dot B^{\frac N2-1}_{2,1}}\Big)\langle t\rangle^{-\frac N4-\frac s2}
\quad  \hbox{ if } \ \ -N/2<s\leq \min\{2,N/2 \},\cr
\|\Lambda^{s}u^{-}\|_{L^2}
\lesssim \Big(D_0+\big\|\big(\nabla R^{+}_{0},\,u^{+}_{0},\,\nabla R^{-}_{0},\,u^{-}_{0}\big)\big\|^h_{\dot B^{\frac N2-1}_{2,1}}\Big)\langle t\rangle^{-\frac N4-\frac s2}
\quad \hbox{ if } \ \ -N/2<s\leq \min\{2,N/2-1 \}.}
$$
This completes the proof of the corollary.
\begin{Remark}\label{1.2} In Corollary \ref{cor1.1},  taking $s=0,1$ leads back to  the standard optimal   $L^{1}$-$L^{2}$ time decay rates \eqref{eq:decate6}-\eqref{eq:decate7} in \cite{EWW} and \eqref{eq:decate10}-\eqref{eq:decate11} in \cite{CWYZ}  when $N=3$,  which is also consistent with the optimal time decay rates \eqref{eq:decate4}-\eqref{eq:decate5}  for a single phase flow model in \cite{Tan1,Tan2}. Note however that our  estimates hold in  the $L^2$  critical framework.
Meanwhile, we also obtain the optimal time  decay rates
 for $(R^{+}-1,\,R^{-}-1)$ as in \eqref{eq:decate1},  which  improves the decay results \eqref{eq:decate8}-\eqref{eq:decate9} in \cite{CWYZ}.
 \end{Remark}

As a consequence
of the   Gagliardo-Nirenberg type inequalities and  Corollary \ref{cor1.1}, similar to the optimal time decay rates \eqref{eq:decate2} and \eqref{eq:decate51} for  single phase flow models,
one can also get the following  more $L^{q}$-$L^{r}$ decay estimates:
\begin{Corollary}\label{cor1.2}

Let the assumptions of Theorem \ref{th:decay} be fulfilled. Then the corresponding solution $(R^{+}-1,\,u^{+},\,R^{-}-1,\,u^{-})$  constructed in  Theorem \ref{th:main1}   satisfies
\begin{eqnarray}
\big\|\Lambda^{k}\big(R^{+}-1,\,u^{+},\,R^{-}-1,\,u^{-}\big)\big\|_{L^p}\lesssim \Big(D_{0}
+\big\|\big(\nabla R^{+}_{0},\,u^{+}_{0},\,\nabla R^{-}_{0},\,u^{-}_{0}\big)\big\|^h_{\dot B^{\frac N2-1}_{2,1}}\Big)\bigl \langle t\rangle^{-\frac N2(1-\frac1p)-\frac k2}, \label{R-E13}
\end{eqnarray}
for all $2\leq p\leq\infty$ and $k\in \mathbb{R}$ satisfying
$-\frac N2<k+N\big(\frac12-\frac1p\big) < \min\big(2,\frac N2-1\big)\cdotp$
\end{Corollary}

\noindent{\bf Proof.} For the following Gagliardo-Nirenberg type inequalities \cite{SS}
\begin{eqnarray*}
\|\Lambda^{k}f\|_{L^{p}}\lesssim \|\Lambda^{m}f\|^{1-\theta}_{L^{q}}\|\Lambda^{n}f\|^{\theta}_{L^{q}},
\end{eqnarray*}
whenever  $0\leq\theta\leq1$,  $1\leq q\leq p\leq\infty$ and $$\ell+N\Big(\frac{1}{q}-\frac{1}{p}\Big)=m(1-\theta)+n\theta,$$
we take $q=2$, $m=\min(2,\frac N2-1)$,
$n=-\frac N2+\epsilon$ with  $\epsilon$ small enough and define $\theta$ by the relation
$$
n\theta+m(1-\theta)=k+N\Bigl(\frac12-\frac1p\Bigr).
$$
Thus,  from  Corollary \ref{cor1.1} we have
\begin{equation*}
\begin{split}
&\big\|\Lambda^{\ell}\big(R^{+}-1,\,u^{+},\,R^{-}-1,\,u^{-}\big)\big\|_{L^{p}}\\
&\quad\lesssim \big\|\Lambda^{m}\big(R^{+}-1,\,u^{+},\,R^{-}-1,\,u^{-}\big)\big\|_{L^2}^{1-\theta}\big\|\Lambda^{n}\big(R^{+}-1,\,u^{+},\,R^{-}-1,\,u^{-}\big)\big\|^{\theta}_{L^2}
\\&\quad \lesssim  \Big(D_{0}
+\big\|\big(\nabla R^{+}_{0},\,u^{+}_{0},\,\nabla R^{-}_{0},\,u^{-}_{0}\big)\big\|^h_{\dot B^{\frac N2-1}_{2,1}}\Big)\Big\{\langle t\rangle^{-\frac N4-\frac m2}\Big\}^{1-\theta}
\Big\{\langle t\rangle^{-\frac N4-\frac n2}\Big\}^{\theta}
\\&\quad= \Big(D_{0}
+\big\|\big(\nabla R^{+}_{0},\,u^{+}_{0},\,\nabla R^{-}_{0},\,u^{-}_{0}\big)\big\|^h_{\dot B^{\frac N2-1}_{2,1}}\Big)\langle t\rangle^{-\frac {N}{4}-\frac{m}{2}(1-\theta)-\frac{n}{2}\theta},
\end{split}
\end{equation*}
 which completes the proof of the corollary.
\begin{center}

\end{center}
\end{document}